\documentclass[12pt,dvipdfmx]{scrartcl}

\usepackage{array}
\usepackage{supertabular}
\usepackage{amsmath,amsfonts,amssymb,amsthm} 
\usepackage{bm}
\usepackage{graphicx}
\usepackage{subcaption}

\newcommand{\nui}{n_i}
\newcommand{\nuip}{n_{i+1}}
\newcommand{\ri}{r_i}
\newcommand{\rip}{r_{i+1}}
\newcommand{\bari}{\bar{r}_i}
\newcommand{\barip}{\bar{r}_{i+1}}
\newcommand{\uni}{U_{(n_i)}}
\newcommand{\unip}{U_{(n_{i+1})}}
\newcommand{\deli}{\Delta_i}
\newcommand{\delip}{\Delta_{i+1}}
\newcommand{\lami}{\lambda_i}
\newcommand{\lamip}{\lambda_{i+1}}

\newcommand{\xwai}{\sum_{i=0}^p}

\allowdisplaybreaks[4]

\newtheorem{theo}{Theorem}
\newtheorem{prop}{Proposition}
\newtheorem{coro}{Corollary}

\title {Estimation of a Continuous Distribution on a Real Line by Discretization Methods--Complete Version--}
\author{Yo Sheena\thanks{Faculty of Economics, Shinshu University; Faculty of Data Science, Shiga University}}

\date{Nov. 2017}

\begin{document}
\maketitle

\begin{abstract}
For an unknown continuous distribution on a real line,  we consider the approximate estimation by the discretization. There are two methods for the discretization. First method is to divide the real line into several intervals before taking samples ("fixed interval method") . Second method  is dividing the real line using the estimated percentiles after taking samples ("moving interval method"). In either way, we settle down to the estimation problem of a multinomial distribution. We use (symmetrized) $f$-divergence in order to measure the discrepancy of the true distribution and the estimated one. Our main result is the asymptotic expansion of the risk (i.e. expected divergence)  up to the second-order term in the sample size. We prove theoretically that the moving interval method is asymptotically superior to the fixed interval method.  We also observe how the presupposed intervals (fixed interval method) or percentiles (moving interval method) affect the asymptotic risk.
\end{abstract}
\noindent
MSC(2010) \textit{Subject Classification}: Primary 60F99; Secondary 62F12
\section{Introduction}
\label{section:int}

One of the useful methods dealing with a continuous distribution is the discretization of the continuous distribution, namely the approximation by the finite-dimensional discrete distribution. Consider a probability distribution on the real line that is absolutely continuous with respect to Lebesgue measure. We call this distribution "mother distribution".  It is not necessarily required to have full support $(-\infty, \infty)$. Let $P(a,b)$ denote the probability of the mother distribution for the interval $(a, b)$. We descretize the mother distribution and get the corresponding multinomial distribution as follows; Let
\begin{equation}
\label{intervals}
-\infty(\triangleq a_0)< a_1 < a_2 < \ldots < a_p < \infty(\triangleq a_{p+1}).
\end{equation}
Consider the multinomial distribution with possible results $C_i\ (i=0,\ldots,p)$ each of which has a probability $P(a_i, a_{i+1})$. This multinomial distribution is an approximation of the mother distribution and coveys a certain amount of information on the mother distribution. In many practical cases, this information could be enough for a statistical analysis with an appropriate selection of $a_i$'s. (See e.g. Drezner and Zerom \cite{Drezner&Zerom} and the cited paper therein for this approximation. )  

In this paper, we consider the estimation of the unknown mother distribution  through this  approximation. Needless to say, the discretized model has a finite number of parameters and much easier to be estimated than the infinite dimensional model for the mother distribution.

There are two methods on how to decide $a_i$'s. One is the "fixed interval method".The $a_i$'s  are given before collecting the sample. In other words, we choose the intervals independently of the sample from the mother distribution.  The other method is the "moving interval method". First choose the percentiles to be estimated $\xi_1 < \ldots < \xi_p$ and estimate them from the sample of the mother distribution. The estimated percentiles $\hat{\xi}_i (i=1,\ldots,p)$ are used as the end points of the intervals, that is, $a_i=\hat{\xi}_i (i=1, \ldots,p)$. The difference between the two methods lies  "intervals first" or "percentiles first".

Once the intervals $a_i$'s are given, we have the estimation problem of the parameters of the multinomial distribution.   If we use the fixed interval method, the true (unknown)  parameters are $P(a_i, a_{i+1}) (i=0, \ldots, p)$ and we need to estimate these parameters based on the sample. On the other hand, for the moving interval method, the true parameter is $P(\hat{\xi}_i, \hat{\xi}_{i+1})$ ($\hat{\xi}_0\triangleq -\infty$ and $\hat{\xi}_{p+1}\triangleq \infty$), while the estimand is the probability given by the presupposed percentiles; if $\xi_i$ is the lower $100\lambda_i\%$ percentile for $1\leq i \leq p$, then the estimated probability for each result is given by $\lambda_{i+1}-\lambda_i (i=0,\ldots,p)$ with $\lambda_0\triangleq 0,\ \lambda_{p+1}\triangleq 1$.

For the measurement of the performance of the estimators, we use $f$-divergence. $f$-divergence between the two multinomial distributions (say $M_1$ and $M_2$) is defined as
\begin{equation}
\label{def_fdive}
D_f[M_1 : M_2] \triangleq
\sum_{i=0}^p p1_i \:f\biggl(\frac{p2_i}{p1_i}\biggr),
\end{equation}
where $p1_i, p2_i, i=0,\ldots,p$ are the probabilities of each result respectively for $M_1$ and $M_2$, and $f$ is a smooth convex  function such that $f(1)=0,\ f'(1)=0,\ f''(1)=1$.  $f$-divergence is natural in view of the sufficiency of the sample information.  If we use the dual function of $f$ defined by $f^*(x)=xf(1/x)$, we have
\begin{equation}
\label{D_f^*}
D_{f^*}[M_1 : M_2] = D_f[M_2 : M_1].
\end{equation}
(See Amari \cite{Amari4} and Vajda \cite{Vajda} for the property of $f$-divergence.)

When the $f$-divergence is too abstract for us to gain some concrete result, we use $\alpha$-divergence.  It is a one-parameter ($\alpha$) family given by \eqref{def_fdive} with $f_\alpha(x)$ such as
\begin{equation}
\label{def_alphadive}
f_\alpha(x)\triangleq 
\begin{cases}
\frac{4}{1-\alpha^2}\bigl(1-x^{(1+\alpha)/2}\bigr)+\frac{2}{1-\alpha}(x-1) &\text{ if $\alpha \ne \pm 1$,} \\
x\log x +1-x & \text{ if $\alpha=1$,}\\
-\log x+x-1 & \text{ if $\alpha=-1$.}
\end{cases}
\end{equation}
We will use the notation $\overset{\alpha}{D}[M_1 : M_2]$ instead of $D_{f_\alpha}[M_1 : M_2]$.
$\alpha$-divergence is the subclass of $f$-divergence, but still a broad class which contains the frequently used divergence such as  Kullback-Leibler divergence ($\alpha =-1$), Hellinger distance ($\alpha=0$), $\chi^2$-divergence ($\alpha=3$). Note that the conjugate of $(f_\alpha)^*$ equals $f_{-\alpha}$, hence 
\begin{equation}
\overset{-\alpha}{D}[M_1 : M_2] = \overset{\alpha}{D}[M_2 : M_1]
\end{equation}

In general, divergence $D[M_1:M_2]$ satisfies the condition
\begin{equation}
D[M_1:M_2]\geq 0,\qquad D[M_1:M_2]=0 \text{ if and only if $M_1\stackrel{d}{=}M_2$}
\end{equation}
But the triangle inequality and symmetricity do not hold true.  In this paper, we adopt the mean of the dual divergences in order to satisfy the symmetricity (see Amari and Cichocki \cite{Amari&Cichocki});
\begin{equation}
\overset{|\alpha|}{D}[M_1 : M_2]\triangleq \frac{1}{2}\left\{
\overset{\alpha}{D}[M_1 : M_2]+\overset{-\alpha}{D}[M_1 : M_2]
\right\}
\end{equation}

We take the expectation of the divergence between the estimated multinomial distribution $\hat{M}$ and the true one $M$;
\begin{equation}
\label{def_risk}
ED \triangleq E\bigl[D_f[M : \hat{M}]\bigr]
\end{equation}
This is the risk of $\hat{M}$ and we use it to describe the goodness of the estimation. In this paper, we only consider the basic estimators, that is, the most likelihood estimator for the fixed interval and the ordered sample for the moving interval.

It is not easy to analyze  the risk theoretically under small sample, hence we focus ourselves on the asymptotic risk under large sample. In Section \ref{section:main}, as the main result, we show the asymptotic expansion of the risk for the both methods, the fixed interval and the moving interval (Theorem 1 and 2 ) . Using this result, first we observe how the asymptotic risk is affected by the the presupposed intervals (the fixed intervals) or percentiles (the moving intervals). Second  we compare the asymptotic risk between the two methods and report the superiority of the moving interval methods when the percentiles are given with equi-probable intervals.
\section{Main Result}
\label{section:main}
We state the asymptotic expansion of the risk \eqref{def_risk} up to the second order with respect to the sample size, $n$, for the both methods, that is, the fixed interval method (Section \ref{subsec:fixed}) and the moving interval method (Section \ref{subsec:moving}). In each subsection, we analyze  how the asymptotic risk is determined with respect to the sample size, the dimension of the multinomial distribution and the prefixed intervals (fixed intervals) or percentiles (moving intervals). In Section \ref{subsec:comparison}, we compare the both methods and show the superiority of the moving intervals when the percentiles are given with equi-probable intervals.
\subsection{Fixed Intervals}
\label{subsec:fixed}
We prefix the intervals with the endpoints \eqref{intervals} \textit{before} taking the sample from the mother distribution. In other words, we choose the endpoints \eqref{intervals} independently of the sample.

We consider the multinomial distribution with the possible results $C_i, i=0,\ldots,p$. If a sample from the mother distribution take the value within the interval $(a_i,  a_{i+1})$ for $i=0,\ldots,p$, we count it as the sample with the result $C_i$. Then this multinomial distribution is an approximation of the mother distribution by a discretizaion..  The probability for $C_i$ is given by
\[
m_i \triangleq P(a_i, a_{i+1}),\quad i=0,\ldots,p,
\]
where $P(a_i, a_{i+1})$ is the probability of the mother distribution for the interval $(a_i, a_{i+1})$.

We estimate this multinomial distribution through the m.l.e.. Let $X_i, i=1,\ldots,n$ be the i.i.d. sample from the mother distribution. Then the m.l.e. of $m\triangleq (m_0, \ldots, m_p)$ is given by $\hat{m}\triangleq (\hat{m}_0, \ldots, \hat{m}_p)$, where 
\begin{equation}
\hat{m}_i\triangleq\# \{X_i | X_i \in (a_i, a_{i+1})\}/n,\quad i=0,\ldots,p.
\end{equation}{dive}, that is, 
\begin{equation}
\label{f-div_fix}
D_f[m: \hat{m}] \triangleq
\sum_{i=0}^p m_i \:f\biggl(\frac{\hat{m_i}}{m_i}\biggr).
\end{equation}
The performance of $\hat{m}$ is measured by the risk, 
\begin{equation}
\label{def_ED_I}
ED_I \triangleq E\bigl[D_f[m : \hat{m}]\bigr].
\end{equation} 

For a general multinomial distribution, which is not necessarily given by a mother distribution as above,  the following result holds.
\begin{theo}
\label{theo_ED_I}
For a multinomial distribution with the probability $m\triangleq(m_0, \ldots, m_p)$ and its m.l.e. $\hat{m}$, the risk of m.l.e. \eqref{def_ED_I}  based on i.i.d. sample of size $n$ is given as follows;
\begin{equation}
\label{ED_I_expan}
ED_I=\frac{p}{2n}+\frac{1}{24n^2}\Bigl[ 4f^{(3)}(1)\Bigl(-3p-1+M \Bigr)+3f^{(4)}(1)\Bigl(-2p-1+M\Bigr)\Bigr],
\end{equation}
where $f^{(3)}$ and $f^{(4)}$ are respectively the third and forth derivative of $f$ in \eqref{f-div_fix},and 
$$
M\triangleq \sum_{i=0}^p m_i^{-1}.
$$
\end{theo}
\noindent
--\textit{Proof}--\\
Let
\[
R_i \triangleq \frac{\hat{m}_i-m_i}{m_i}.
\]
Note that 
\[
\Bigl(\sqrt{n}(\hat{m}_1-m_1),\ldots, \sqrt{n}(\hat{m}_p-m_p)\Bigr) \stackrel{d}{\longrightarrow} N_p(0, \Sigma),
\]
where
\[
\Sigma\triangleq (\sigma_{ij}),\qquad \sigma_{ij}\triangleq
\begin{cases}
p_i(1-p_i) &\text{ if $i=j$, }\\
-p_i p_j&\text{ if $i \ne j$.}
\end{cases}
\]
(See e.g. (5.4.13) of \cite{Lehmann}.)  Using this fact and $f(1)=0,\ f'(1)=0,\ f''(1)=1$, we have the following expansion $D_f[m: \hat{m}]$ with respect to $n$. 
\begin{align}
\label{expan_D_f}
&D_f[m: \hat{m}]\nonumber\\
&=\sum_{i=0}^p m_i f(1+R_i)\nonumber\\
&=\sum_{i=0}^p m_i \Bigl(f(1)+f'(1)R_i+\frac{1}{2}f''(1)R_i^2+\frac{1}{6}f^{(3)}(1)R_i^3+\frac{1}{24}f^{(4)}(1)R_i^4\Bigr)+o_p(n^{-2})\nonumber\\
&=\frac{1}{2}\sum_{i=0}^p m_i R_i^2 +\frac{1}{6}f^{(3)}(1)\sum_{i=0}^p m_i R_i^3+\frac{1}{24}f^{(4)}(1)\sum_{i=0}^p m_i R_i^4+o_p(n^{-2}).\nonumber\\
&=\frac{1}{2}\sum_{i=0}^p m_i^{-1}(\hat{m}_i-m_i)^2  +\frac{1}{6}f^{(3)}(1)\sum_{i=0}^p m_i^{-2}(\hat{m}_i-m_i)^3\nonumber\\
&\qquad+\frac{1}{24}f^{(4)}(1)\sum_{i=0}^p m_i^{-3}(\hat{m}_i-m_i)^4+o_p(n^{-2}).
\end{align}
From the central moments of the standardized multinomial distribution,
\begin{align*}
&E[\hat{m}_i-m_i]=0,\qquad E[(\hat{m}_i-m_i)^2]=n^{-1}(m_i-m_i^2), \\
&E[(\hat{m}_i-m_i)^3]=n^{-2}(m_i-3m_i^2+2m_i^3),\qquad E[(\hat{m}_i-m_i)^4]=3n^{-2}(m_i-m_i^2)^2+o(n^{-2}),
\end{align*}
we have
\[
ED_I=\frac{1}{2n}\sum_{i=0}^p(1-m_i)+\frac{1}{6n^2}f^{(3)}(1)\sum_{i=0}^p(m_i^{-1}-3+2m_i)+\frac{1}{8n^2}f^{(4)}(1)\sum_{i=0}^p(m_i^{-1}-2+m_i),
\]
which is equivalent to the result \eqref{ED_I_expan} since $\sum_{i=0}^p m_i=1$.\hfill \textit{ Q.E.D.}

Especially for the $\alpha$-divergence,
\[
\overset{\alpha}{D}[m : \hat{m}]\triangleq D_{f_\alpha}[m : \hat{m}] ,\qquad \overset{|\alpha|}{D}[m : \hat{m}]\triangleq \frac{1}{2}\Bigl\{ D_{f_\alpha}[m : \hat{m}]+D_{f_{-\alpha}}[m : \hat{m}]\Bigr\},
\]
where $f_\alpha$ is given by \eqref{def_alphadive}, 
the following results hold. (Sheena \cite{Sheena} gained this result as an example of  the asymptotic risk of m.l.e. for a general parametric model.)
\begin{coro}
\begin{align}
\label{ED_I_alpha_expan}
\overset{\alpha}{ED}_I&\triangleq E\bigl[\overset{\alpha}{D}[m: \hat{m}]\bigr]=\frac{p}{2n}+\frac{1}{96n^2}\Bigl\{(\alpha-3)(3\alpha-7)(M-1)-6(\alpha-3)(\alpha-1)p\Bigr\}+o(n^{-2}),\\
\label{ED_I_alpha_abso_expan}
\overset{|\alpha|}{ED}_I&\triangleq E\bigl[\overset{|\alpha|}{D}[m: \hat{m}]\bigr]=\frac{p}{2n}+\frac{1}{32n^2}\Bigl\{(\alpha^2+7)(M-1)-2(\alpha^2+3)p\Bigr\}+o(n^{-2}).
\end{align}
\end{coro}
\noindent --\textit{Proof}--
\\
The results are straightforward from Theorem \ref{theo_ED_I} and the fact 
\begin{equation}
\label{f_deri}
f^{(3)}_\alpha(1)=(\alpha-3)/2 \qquad f^{(4)}_\alpha(1)=(\alpha-3)(\alpha-5)/4.
\end{equation}
\hfill \textit{ Q.E.D.}

We observe the following points from \eqref{ED_I_expan}, \eqref{ED_I_alpha_expan} and \eqref{ED_I_alpha_abso_expan}.
\begin{enumerate}
\item The main term, i.e. $n^{-1}$-order term, is determined by $p/n$, that is the ratio of the dimension of the multinomial distribution model (the number of the free parameters) to the sample size. We call this  "$p-n$ ratio" hereafter. $p-n$ ratio shows  the  complexity of the model to be estimated relative to the  sample size. The main term is independent of  $f$ or $\alpha$, and $m_i (i=0, \ldots, p)$.
\item The second term, i.e. $n^{-2}$-order term, depends on the parameter of the multinomial distribution through 
\[
M\triangleq\sum_{i=0}^{p} m_i^{-1}.
\]
$M$ attains the minimum value $(p+1)^2$ when $m_0=m_1=\cdots=m_p$. It increases rapidly if one of  $m_i$'s is near to zero. The effect of $M$ on the risk depends on the choice of $f$ or $\alpha$. If you choose $f$ such that $4f^{(3)}(1)+3f^{(4)}(1)$ is non-positive or $\alpha$ such that $7/3 \leq \alpha \leq 3$, \eqref{ED_I_expan} and \eqref{ED_I_alpha_expan} respectively decreases or are constant as $M$ increases. This is rather unnatural since it contradicts to our \textit{belief} that the existence of result with a small probability makes estimation harder for a multinomial distribution. In this sense, $\chi^2$-distance with $\alpha=3$ seems inappropriate, since it is asymptotically insensitive to the difference in the parameters $m_i (i=0, \ldots, p)$.  (See Sheena \cite{Sheena_2}, which reports that the $\alpha$-divergence seems statistically unnatural when $|\alpha|$ is large for a regression model.) $\alpha$-divergence is a distance if and only if $\alpha=0$, and the pair of $\alpha$- and $-\alpha$- divergences work dually like a distance. (For "generalized Pythagorean theorem", see \cite{Amari4} or \cite{Amari&Nagaoka}.) In this respect, the divergence $\overset{|\alpha|}{D}$ seems natural. Actually \eqref{ED_I_alpha_abso_expan} shows that the risk is a monotonically increasing function of $M$ for any $\alpha$. 
\item The $n^{-2}$ term of \eqref{ED_I_expan} or \eqref{ED_I_alpha_expan} can be negative for some $f(\text{or $\alpha$}), p$ , while that of \eqref{ED_I_alpha_abso_expan} is always positive as
\[
(\alpha^2+7)(M-1)-2(\alpha^2+3)p \geq  (\alpha^2+7)((p+1)^2-1)-2(\alpha^2+3)p=p^2\alpha^2+7p^2+8p>0.
\]
\end{enumerate}
\subsection{Moving Intervals}
\label{subsec:moving}
First we choose points $\lambda_i (1\leq i \leq p)$ in the interval $(0, 1)$;
\begin{equation}
\label{def_per}
\lambda_0(\triangleq 0)< \lambda_1 < \lambda_2 < \cdots < \lambda_p <\lambda_{p+1}(\triangleq 1).
\end{equation}
Let 
\begin{equation}
\label{def_per_mother}
\xi_i \triangleq F^{-1}(\lambda_i),\ 1\leq i \leq p, \quad \xi_0\equiv -\infty,\quad \xi_{p+1}\equiv \infty,
\end{equation}
where $F^{-1}$ is the inverse function of the cumulative distribution function, $F$, of the mother distribution. 
We call $\xi$'s the percentiles of the mother distribution. 

In the moving intervals method, we estimate the percentiles of the mother distribution from the sample of the mother distribution, and use them as the endpoints of \eqref{ED_I_alpha_abso_expan};
\begin{equation}
a_i = \hat{\xi}_i,\quad 1\leq i \leq p,
\end{equation}
where $\hat{\xi}_i$ is the estimator of $\xi_i$ for $i=1,\ldots p$ and $\hat{\xi}_0\equiv -\infty$ and $\hat{\xi}_{p+1}\equiv \infty$.  In this case,  the multinomial distribution that approximates the mother distribution has unknown parameters 
\[
\hat{m}\triangleq(\hat{m}_0,\ldots,\hat{m}_p),\qquad \hat{m}_i\triangleq P(a_{i}, a_{i+1}) \equiv P(\hat{\xi}_i, \hat{\xi}_{i+1}) \quad 0\leq i \leq p,
\]
while it is estimated as 
\begin{equation}
\label{def_m_moving}
m\triangleq(m_0, \ldots, m_p),\qquad m_i\triangleq \lambda_{i+1}-\lambda_{i} \quad 0\leq i \leq p.
\end{equation}

Although there are several ways to estimate the percentile $\xi$, we focus here on the simple estimator using the order statistic itself.
Take i.i.d sample of size $n$ from the mother distribution and let the ordered sample be denoted by
\[
X_{(1)} \leq X_{(2)} \leq \cdots \leq X_{(n)}.
\]
We estimate $\xi _i $ by 
\begin{equation}
\hat{\xi}_i\triangleq X_{(n_i)}\quad 1\leq i \leq p,
\end{equation}
where $n_i$ is a function of $n$ with the values in $\{1, 2, \ldots, n\}.$  
Let $r_i$ denote the gap between $n_i$ and $n\lambda_i$, namely
\begin{equation}
\label{def_ri}
r_i\triangleq n_i-n\lambda_i \quad 1\leq i \leq p, \qquad r_0\triangleq 0, \qquad r_{p+1}\triangleq 1.
\end{equation}

We measure the discrepancy between $m$ and $\hat{m}$ by $f$-divergence, 
\begin{equation}
\label{f-div_move}
D_f[m : \hat{m}] \triangleq
\sum_{i=0}^p m_i \:f\biggl(\frac{\hat{m}_i}{m_i}\biggr).
\end{equation}
If one might think it is natural to consider $D_f[\hat{m} : m]$ in the sense that the true parameter should come first, it is satisfied by using the dual function $f^*$ (see \eqref{D_f^*}). Hence we will  proceed with \eqref{f-div_move}.

The risk for the moving interval method is given by
\begin{equation}
\label{def_ED_P}
ED_P \triangleq E\bigl[D_f[m : \hat{m}]\bigr],
\end{equation} 
and the following result holds.
\begin{theo}
Suppose that $r_i/n=o(n_i^{-1/2})$, then 
\begin{align}
ED_P&=\frac{p}{2n}+\frac{1}{24 n^2}\Bigl[-24-36p+12\sum_{i=0}^p(r_{i+1}-r_i)(r_{i+1}-r_i+1)m_i^{-1}\nonumber\\
&\hspace{30mm}+4f^{(3)}(1)\Bigl\{-5-9p+\sum_{i=0}^p\bigl(3(r_{i+1}-r_i)+2\bigr)m_i^{-1}\Bigr\}\nonumber\\
&\hspace{30mm}+f^{(4)}(1)\Bigl\{-3-6p+3\sum_{i=0}^p m_i^{-1}\Bigr\}\Bigr]+o(n^{-2}).\label{ED_P_expan}
\end{align}
\end{theo}
\noindent --\textit{Proof}--
\\
The whole process of proof is lengthy, hence we only state the outline of the proof here. All the details are found in Appendix. 
Let
\[
U_{(n_i)}\triangleq F(X_{(n_i)}), \qquad \Delta_i \triangleq \sqrt{n}(U_{(n_i)}-\lambda_i) \quad 1\leq i \leq p
\]
and $\Delta_0\triangleq 0,\ \Delta_{p+1}\triangleq 0$. The following relationship holds for $0 \leq i \leq p$.
\begin{align}
\hat{m}_i&=F(\hat{\xi}_{i+1})-F(\hat{\xi}_i)\nonumber\\
&=F(X_{(n_{i+1})})-F(X_{(n_i)})\nonumber\\
&=U_{(n_{i+1})}-U_{(n_i)}\nonumber\\
&=\lambda_{i+1}-\lambda_i+n^{-1/2}(\Delta_{i+1}-\Delta_i)\nonumber\\
&=m_i+n^{-1/2}(\Delta_{i+1}-\Delta_i). \label{form_hat_m}
\end{align}
Note that
\[
(\Delta_1, \ldots, \Delta_p) \stackrel{d}{\longrightarrow} N_p(0, \Sigma),
\]
where 
\[
\Sigma=(\sigma_{ij})=\lambda_i (1-\lambda_j)\quad 1\leq i \leq j \leq p
\]
(see e.g. Theorem 5.4.5 of \cite{Lehmann}). Similarly to \eqref{expan_D_f}, the following equation holds. 
\begin{equation}
D_f[m: \hat{m}]=\frac{1}{2}\sum_{i=0}^p m_i R_i^2 +\frac{1}{6}f^{(3)}(1)\sum_{i=0}^p m_i R_i^3+\frac{1}{24}f^{(4)}(1)\sum_{i=0}^p m_i R_i^4+o_p(n^{-2}).
\end{equation}
Therefore we have
\begin{equation}
\label{ED_P_pro}
ED_P=\frac{1}{2}\sum_{i=0}^p m_i E[R_i^2] +\frac{1}{6}f^{(3)}(1)\sum_{i=0}^p m_i E[R_i^3]+\frac{1}{24}f^{(4)}(1)\sum_{i=0}^p m_i E[R_i^4]+o(n^{-2}).
\end{equation}
After long but straightforward calculation (see Appendix), we have
\begin{align}
\sum_{i=0}^p m_i E[R_i^2]&=n^{-1}p+n^{-2}[-2-3p+\sum_{i=0}^p (r_{i+1}-r_i)(r_{i+1}-r_i+1)m_i^{-1}],\label{E_R2}\\
\sum_{i=0}^p m_i E[R_i^3]&=n^{-2}[-5-9p+\sum_{i=0}^p\bigl(3(r_{i+1}-r_i)+2\bigr)m_i^{-1}],\label{E_R3}\\
\sum_{i=0}^p m_i E[R_i^4]&=n^{-2}[-3-6p+3\sum_{i=0}^pm_i^{-1}].\label{E_R4}
\end{align}
If we insert these results into \eqref{ED_P_pro}, we have the result.
\hfill \textit{ Q.E.D.}

We also have the following formulas for the $\alpha$-divergence.
\begin{coro}
\begin{align}
\label{ED_P_alpha_expan}
\overset{\alpha}{ED}_P&=\frac{p}{2n}+\frac{1}{96n^2}\Bigl[-\alpha^2(3+6p)-\alpha(16+24p)-18p-21\nonumber\\
&\hspace{25mm}+\sum_{i=0}^p\bigl\{48(r_{i+1}-r_i)^2+24(\alpha-1)(r_{i+1}-r_i)+3\alpha^2-8\alpha-3\bigr\}m_i^{-1}\Bigr]\nonumber\\
&\quad+o(n^{-2}),\\
\label{ED_P_alpha_abso_expan}
\overset{|\alpha|}{ED}_P&=\frac{p}{2n}+\frac{1}{96n^2}\Bigl[-\alpha^2(3+6p)-18p-21\nonumber\\
&\hspace{25mm}+\sum_{i=0}^p\bigl\{48(r_{i+1}-r_i)^2-24(r_{i+1}-r_i)+3\alpha^2-3\bigr\}m_i^{-1}\Bigr]+o(n^{-2})
\end{align}
\end{coro}
\noindent --\textit{Proof}--
\\
The results are straightforward from \eqref{ED_P_expan} and \eqref{f_deri}.
\hfill \textit{ Q.E.D.}

We give some cements on $ED_P$, $\overset{\alpha}{ED}_P$ and $\overset{|\alpha|}{ED}_P$.
\begin{enumerate}
\item The main term is half the $p-n$ ratio just like $ED_I$. It is independent of  $f$ or $\alpha$, and $m_i (i=0, \ldots, p)$.
\item The risk is independent of the mother distribution (it is due to the fact \eqref{form_hat_m}).  It is determined by our choice of $m_i$'s or equivalently $\lambda_i$'s in \eqref{def_per}. 
\item The choice of $n_i$'s, or equivalently $r_i$'s $(i=1, \ldots, p)$ effects the $n^{-2}$-order term. It is possible that the coefficient of $m_i^{-1}$ could be negative for some $r_i$'s and $f$(or $\alpha$). In this case small $m_i$ could reduce the risk. 
\end{enumerate}
\subsection{Comparison of two methods}
\label{subsec:comparison}
We compare the risks between the fixed interval method and the moving interval method. For the both methods, the main term ($n^{-1}$-order term) are common, but we can see some difference in the second term ($n^{-2}$-order term). The biggest difference between the two methods lies in $m_i$'s.  In the fixed interval method, $m_i$'s depend on the unknown mother distribution, hence we are unable to control them. As we observed in Section \ref{subsec:fixed}, if they include even one small $m_i$ near to zero, then the (asymptotic) risk gets extremely high through $M$.  The  more intervals (endpoints)  we use for discretization, more likely we are to have small $m_i$'s. Even if we have a large set of sample, we have to be cautions to raise the dimension of the multinomial distribution. On the contrary, for the moving interval method, $m_i$'s are controllable. We can choose $m_i$'s so that the risk does not take a large value.

In order to make more specific comparison,  first we will specify $n_i$'s or equivalently $r_i$'s for the moving interval method.
The most naive selection of $n_i$ is $[n\lambda_i]$ or $[n\lambda_i]+1$, where $[\ \cdot \ ]$ is Gauss symbol.
Let
\begin{equation}
\bari\triangleq [n\lambda_i]-n\lambda_i.
\end{equation}
In this paper, we adopt the following randomized choice of $r_i$'s;
\begin{equation}
\label{dist_ri}
P(r_i=\bari)\Bigl(=P(n_i=[n\lambda_i])\Bigr)=1+\bari,\quad P(r_i=1+\bari)\Bigl(=P(n_i=[n\lambda_i]+1)\Bigr)=-\bari
\end{equation}
for $1\leq i \leq p$, while $r_0\equiv 0$ and $r_{p+1}\equiv1$ as in \eqref{def_ri}.
This is natural in that $n_i$ is chosen to be $[n\lambda_i]$ and $[n\lambda_i]+1$ respectively with the probabilities proportional to the closeness to the both points. (To locate $\hat{\xi}_i$ between $X_{([n\lambda_i])}$ and $X_{([n\lambda_i]+1)}$ according to $r_i$ is another appealing idea. But if we adopt this estimation of $\xi_i$, then the risk depends on the mother distribution.)

Let 
\[
ED_P^*\triangleq E[ED_P],\quad \overset{\alpha}{ED^*_P}\triangleq E[\overset{\alpha}{ED}_P],\quad \overset{|\alpha|}{ED_P^*}\triangleq E[\overset{|\alpha|}{ED}_P],
\]
where all the expectation is taken with respect to the distribution \eqref{dist_ri}.
The following results hold for the randomized choice of $r_i$'s \eqref{dist_ri}.
\begin{prop}
\begin{align}
\label{ED_P^*_expan}
ED_P^*&=\frac{p}{2n}+\frac{1}{48n^2}\Bigl[-48-72p+24\Bigl\{-\bar{r}_1(1+\bar{r}_1)m_0^{-1}+\bigl(2-\bar{r}_p(1+\bar{r}_p)\bigr)m_p^{-1}\nonumber\\
&\hspace{60mm}-\sum_{i=1}^{p-1}\bigl(\bari(1+\bari)+\barip(1+\barip)\bigr)m_i^{-1}\Bigr\}\nonumber\\
&\hspace{30mm}+8f^{(3)}(1)\Bigl\{-5-9p+2\sum_{i=0}^p m_i^{-1}+3m_p^{-1}\Bigr\}\nonumber\\
&\hspace{30mm}+2f^{(4)}(1)\Bigl\{-3-6p+3\sum_{i=0}^p m_i^{-1}\Bigr\}\Bigr]+o(n^{-2}),\\
\label{ED_P^*_alpha_expan}
\overset{\alpha}{ED_P^*}&=\frac{p}{2n}+\frac{1}{96n^2}\Bigl[-\alpha^2(3+6p)-\alpha(16+24p)-18p-21\nonumber\\
&\hspace{30mm}-48\bar{r}_1(1+\bar{r}_1)m_0^{-1}+\bigl(-48\bar{r}_p(1+\bar{r}_p)+24(\alpha+1)\bigr)m_p^{-1}\nonumber\\
&\hspace{30mm}-48\sum_{i=1}^{p-1}\bigl(\bari(1+\bari)+\barip(1+\barip)\bigr)m_i^{-1}\nonumber\\
&\hspace{30mm}+(3\alpha^2-8\alpha-3)\sum_{i=0}^p m_i^{-1}\Bigr]+o(n^{-2}),\\
\label{ED_P^*_alpha_abso_expan}
\overset{|\alpha|}{ED_P^*}&=\frac{p}{2n}+\frac{1}{96n^2}\Bigl[-\alpha^2(3+6p)-18p-21\nonumber\\
&\hspace{30mm}-48\bar{r}_1(1+\bar{r}_1)m_0^{-1}+\bigl(24-48\bar{r}_p(1+\bar{r}_p)\bigr)m_p^{-1}\nonumber\\
&\hspace{30mm}-48\sum_{i=1}^{p-1}\bigl(\bari(1+\bari)+\barip(1+\barip)\bigr)m_i^{-1}\nonumber\\
&\hspace{30mm}+(3\alpha^2-3)\sum_{i=0}^p m_i^{-1}\Bigr]+o(n^{-2}).
\end{align}
\end{prop}
\noindent --\textit{Proof}--
\\
As proved in Appendix, the following results hold.
\begin{align}
&E[\ri]=0  \text{ for $0 \leq i \leq p$}, \qquad E[r_{p+1}]=1 \label{exp_ri}\\
&E[\ri^2]=-\bari(1+\bari)  \text{ for $1\leq i \leq p$}, \qquad E[r_0^2]=0, \qquad E[r_{p+1}^2]=1 \label{exp_ri^2}\\
&E[\ri\rip]=0 \text{ for $0\leq i \leq p$}. \label{exp_ri_rip}
\end{align}
Applying these results to $E[(\rip-\ri)^2]=E[\ri^2]+E[\rip^2]-2E[\ri\rip]$ and $E[\rip]-E[\ri]$ in \eqref{ED_P_expan}, \eqref{ED_P_alpha_expan} and \eqref{ED_P_alpha_abso_expan}, we have the results. \hfill \textit{ Q.E.D.}
\\

Note that for $1 \leq i \leq p$, $-1 \leq \bari \leq 0$ and 
$$
 0 \leq -\bari(1+\bari) \leq  1/4. 
$$
Therefore we have
\begin{align}
\label{ED_P^*_upbound}
ED_P^*&\leq \frac{p}{2n}+\frac{1}{48n^2}\Bigl[-48-72p+6\Bigl\{m_0^{-1}+9m_p^{-1}+2\sum_{i=1}^{p-1}m_i^{-1}\Bigr\}\nonumber\\
&\hspace{30mm}+8f^{(3)}(1)\Bigl\{-5-9p+2\sum_{i=0}^p m_i^{-1}+3m_p^{-1}\Bigr\}\nonumber\\
&\hspace{30mm}+2f^{(4)}(1)\Bigl\{-3-6p+3\sum_{i=0}^p m_i^{-1}\Bigr\}\Bigr]+o(n^{-2}) \quad\Bigl(say\ \overline{ED}_P^*\Bigr),\\
\label{ED_P^*_alpha_upbound}
\overset{\alpha}{ED_P^*}&\leq \frac{p}{2n}+\frac{1}{96n^2}\Bigl[-\alpha^2(3+6p)-\alpha(16+24p)-18p-21\nonumber\\
&\hspace{30mm}+12m_0^{-1}+(24\alpha+36)m_p^{-1}\nonumber\\
&\hspace{30mm}+24\sum_{i=1}^{p-1}m_i^{-1}+(3\alpha^2-8\alpha-3)\sum_{i=0}^p m_i^{-1}\Bigr]+o(n^{-2})\nonumber\\
&=\frac{p}{2n}+\frac{1}{96n^2}\Bigl[-\alpha^2(3+6p)-\alpha(16+24p)-18p-21\nonumber\\
&\hspace{30mm}+(3\alpha^2-8\alpha+9)m_0^{-1}+(3\alpha^2+16\alpha+33)m_p^{-1}\nonumber\\
&\hspace{30mm}+(3\alpha^2-8\alpha+21)\sum_{i=1}^{p-1}m_i^{-1}\Bigr]+o(n^{-2})
\quad\Bigl(say\ \overset{\alpha}{\overline{ED}_P^*}\Bigr), \\
\label{ED_P^*_alpha_abso_upbound}
\overset{|\alpha|}{ED_P^*}&\leq\frac{p}{2n}+\frac{1}{32n^2}\Bigl[-\alpha^2(1+2p)-6p-7+(\alpha^2+3)m_0^{-1}+(\alpha^2+11)m_p^{-1}\nonumber\\
&\hspace{30mm}+(\alpha^2+7)\sum_{i=1}^{p-1}m_i^{-1}\Bigr]+o(n^{-2})
\quad\Bigl(say\ \overset{|\alpha|}{\overline{ED}_P^*}\Bigr).
\end{align}
If we choose the equal right-end and left-end probabilities, i.e. $m_0=m_{p}$, 
\begin{equation}
\overset{|\alpha|}{\overline{ED}_P^*}=\frac{p}{2n}+\frac{1}{32n^2}\Bigl[-\alpha^2(1+2p)-6p-7+(\alpha^2+7)M\Bigr]+o(n^{-2}).
\end{equation}
This upper bound for $\overset{|\alpha|}{ED_P^*}$  is affected by $m_i$'s through $M$ just like \eqref{ED_I_expan}. This indicates that the choice of equally-valued $m_i$'s, that is, $m_i=1/(p+1), i=1, \ldots, p$ are reasonable for the estimation of the mother distribution. It is needles to say that the percentiles with a common increment ("quantiles") are most often used in a practical situation.  If we choose "quantiles" for the moving interval method, we have the following result.
\begin{theo}
Set $\lambda_i$'s in \eqref{def_per} so that $m_i=1/(p+1), \ i=0,\ldots,p$, then asymptotically (exactly speaking, as for the comparison up to the $n^{-2}$-order term) , the following inequality holds.
\begin{equation}  
\overset{|\alpha|}{ED}_I \geq \overset{|\alpha|}{ED_P^*}.
\end{equation}
\end{theo}
\noindent --\textit{Proof}--
\\
Since $M \geq (p+1)^2$, from \eqref{ED_I_alpha_abso_expan}, we have
$$
\overset{|\alpha|}{ED}_I  \geq \frac{p}{2n}+\frac{1}{32n^2}\Bigl\{(\alpha^2+7)(p^2+2p)-2(\alpha^2+3)p\Bigr\}+o(n^{-2}) \quad \Bigl( say\  \overset{|\alpha|}{\underline{ED}}_I \Bigr),
$$
while when $m_i=1/(p+1),\ i=0, \ldots, p$, $\overset{|\alpha|}{\overline{ED}_P^*}$ equals
$$
\frac{p}{2n}+\frac{1}{32n^2}\Bigl[-\alpha^2(1+2p)-6p-7+(\alpha^2+7)(p^2+2p+1)\Bigr]+o(n^{-2}).
$$
Up to the $n^{-2}$-order term, we have
\begin{equation}
\overset{|\alpha|}{ED}_I-\overset{|\alpha|}{ED_P^*}\geq \overset{|\alpha|}{\underline{ED}}_I - \overset{|\alpha|}{\overline{ED}_P^*} =0.
\end{equation}
\hfill \textit{ Q.E.D.}
\\

The above theorem says that even if we are lucky enough to choose the best intervals (that is, equi-probable intervals) for the fixed interval method, it is asymptotically dominated by the moving interval method with "quantiles". We can conclude that if we estimate an unknown continuous distribution by the approximation method of discretization, it is better, at least asymptotically, to use the moving interval method.

We will also present a numerical comparison between the both methods. Suppose that $a_i$'s in \eqref{intervals} for the fixed interval method is given by
\begin{equation}
\label{example_a_i}
(-2.0, -1.5, -1.0, -0.5,\ 0,\  0.5, \ 1.0,\ 1.5, \ 2.0).
\end{equation}
with $p=9$. 
We consider the two cases where the mother distribution are respectively $N(0,1)$ and $st(0.8)$, where $st(0.8)$ is the skew $t$-distribution with the zero mean, the unit variance and the skewness parameter of 0.8.

For the intervals with the endpoints \eqref{example_a_i}, the corresponding probabilities of $N(0,1)$ are
$$
(m_0, m_1, \ldots, m_9) \doteqdot (0.023, 0.044, 0.092, 0.150, 0.191, 0.191, 0.150, 0.092, 0.044, 0.023),
$$
while those of $st(0.8)$ are given by
$$
(m_0, m_1, \ldots, m_9) \doteqdot (6.496*10^{-8}, 0.003, 0.153, 0.219, 0.194, 0.155, 0.113, 0.074, 0.044, 0.044)
$$

\begin{figure}
\centering
\begin{subfigure}{0.4\columnwidth}
\centering
\includegraphics[width=\columnwidth]{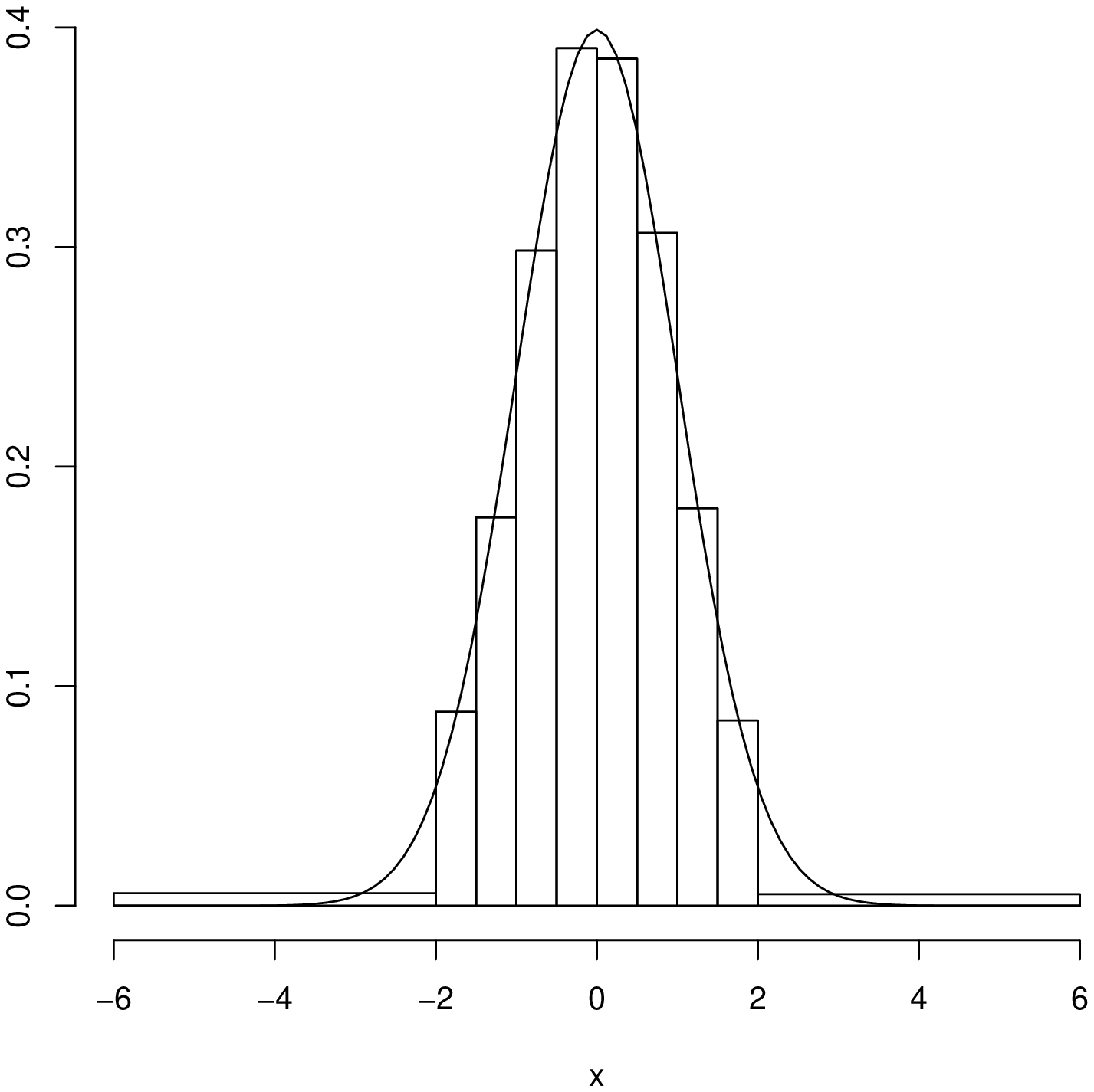}
\caption{$N(0,1)$}
\label{dens_N}
\end{subfigure}
\hspace{5mm}
\begin{subfigure}{0.4\columnwidth}
\centering
\includegraphics[width=\columnwidth]{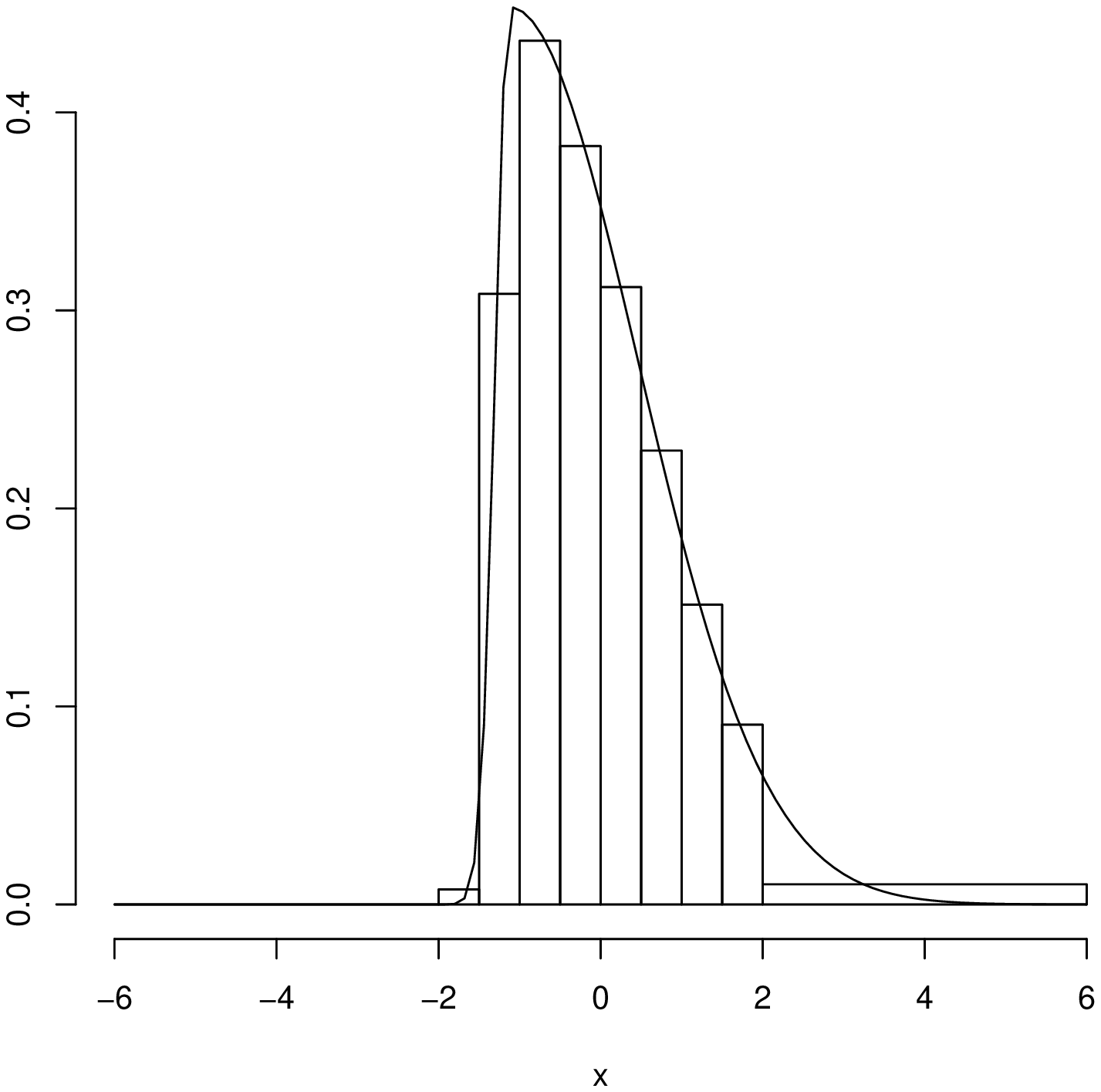}
\caption{$st(0.8)$}
\label{dens_t}
\end{subfigure}
\caption{Density and Histogram of $10^4$ samples}
\label{dens_hist_mother}
\end{figure}

The density function of $N(0,1)$  and the histogram of $10^4$ samples with the above endpoints \eqref{example_a_i} are drawn in Figure \ref{dens_N}. The similar figures for $st(0.8)$ are drawn in Figure \ref{dens_t}.

For the moving interval method, we use "quantiles". Namely  $\lambda$'s in \eqref{def_per} are given by $\lambda_i=i/10 \  (1 \leq i \leq 9)$, or equivalently $m_i (0\leq i \leq 9)$ in \eqref{def_m_moving} are all $1/10$. 

\begin{figure}
\centering
\begin{subfigure}{0.4\columnwidth}
\centering
\includegraphics[width=\columnwidth]{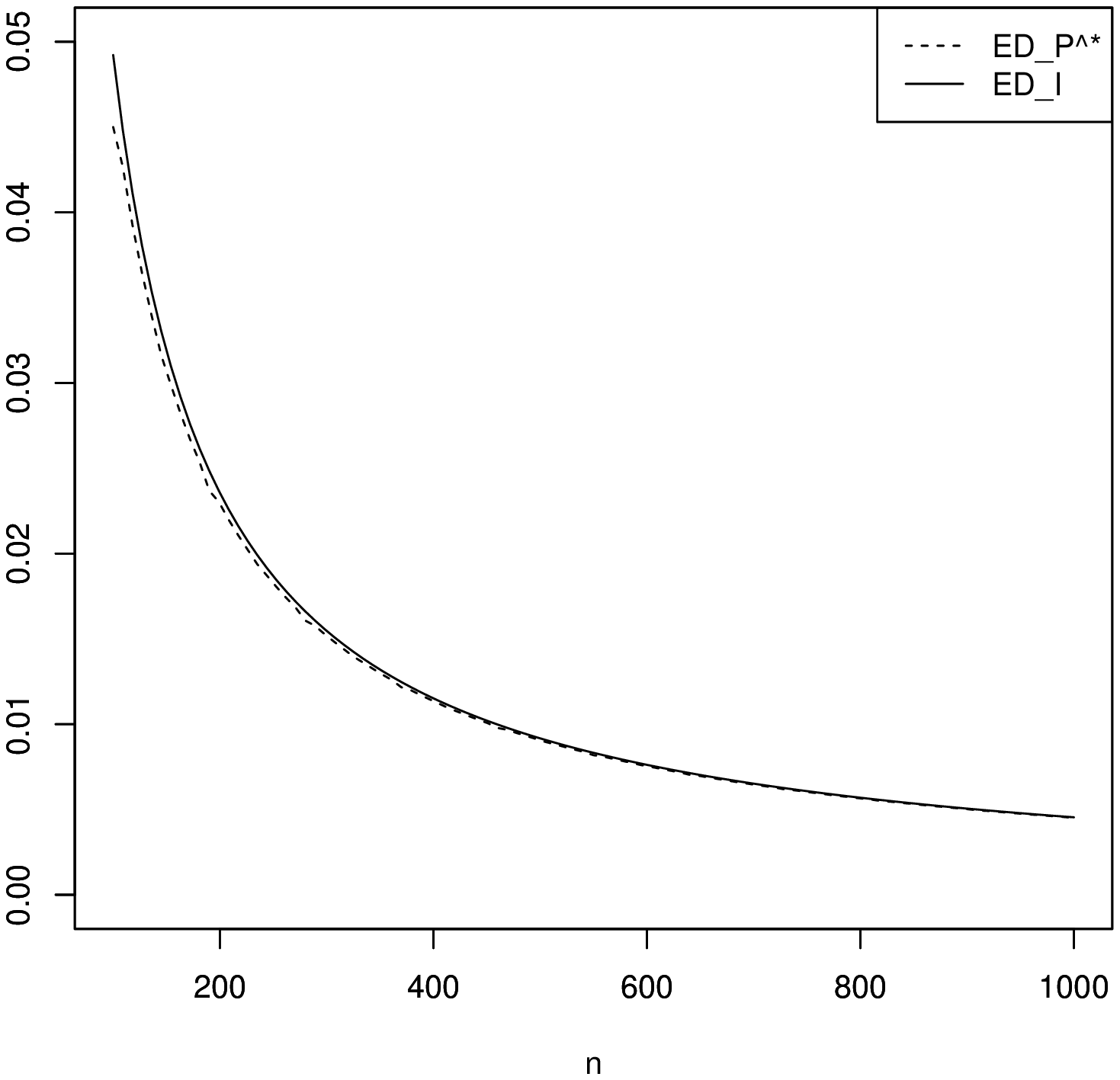}
\caption{N(0,1)}
\label{risks_N}
\end{subfigure}
\hspace{5mm}
\begin{subfigure}{0.4\columnwidth}
\centering
\includegraphics[width=\columnwidth]{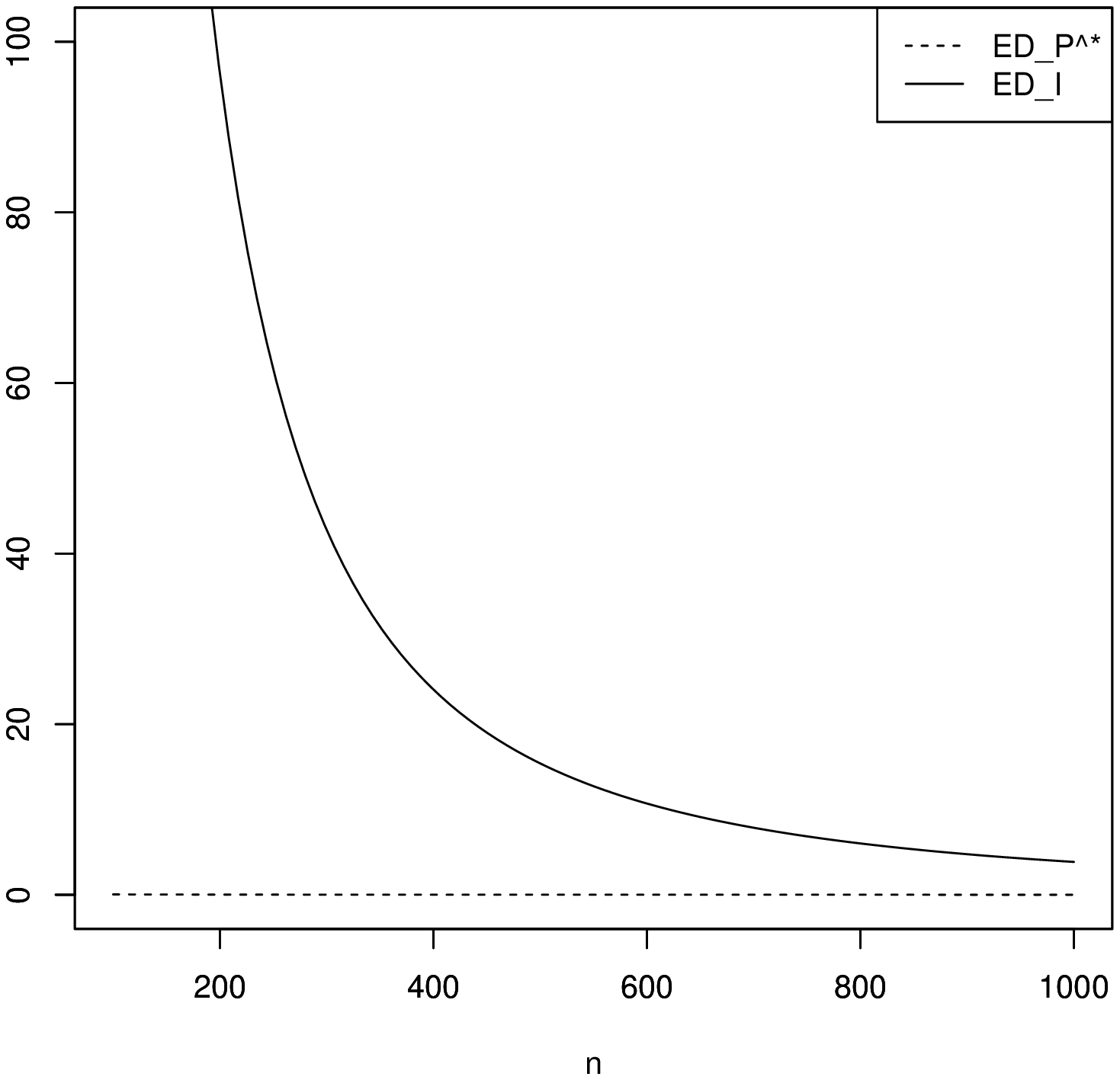}
\caption{t(0.8)}
\label{risks_t}
\end{subfigure}
\caption{The risks for the both methods}
\label{risks_infty_int}
\end{figure}

We put $\alpha=1$. Let's skip the $o(n^{-2})$ part of $\overset{|\alpha|}{ED}_I$ ($\overset{|\alpha|}{ED_P^*}$), and call it the approximated  $\overset{|\alpha|}{ED}_I$ ($\overset{|\alpha|}{ED_P^*}$). The graphs of the approximated $\overset{|\alpha|}{ED}_I$ and $\overset{|\alpha|}{ED_P^*}$ as $n$ varies are drawn in Figure \ref{risks_N} for $N(0,1)$ and in Figure \ref{risks_t} for $st(0.8)$. (Note that $|\alpha|$ are skipped in the legend.)
In Figure \ref{risks_N}, though the graph of the approximated $\overset{|\alpha|}{ED_P^*}$ is slightly lower than that of $\overset{|\alpha|}{ED}_I$, the two curves are quite close to each other. In Figure \ref{risks_t}, we see that the curve of $\overset{|\alpha|}{ED}_I$ is located at much higher position than that of $\overset{|\alpha|}{ED_P^*}$.

Let's  consider the approximated $\overset{|\alpha|}{ED}_I$ and $\overset{|\alpha|}{ED_P^*}$ as the functions of $n$ and put the equation
\begin{equation}
\label{eq_EDI_EDP}
\text{ (The approximated $\overset{|\alpha|}{ED}_I)(n)$}=\text{ (The approximated $\overset{|\alpha|}{ED_P^*})(100)$}
\end{equation}
The solution of this equation indicates how large sample is required for the approximated $\overset{|\alpha|}{ED}_I$ to attain the same risk as that of the approximated $\overset{|\alpha|}{ED_P^*}$ with $n=100$. For the case of $N(0,1)$ the solution is given by $n\doteqdot 109$, while $n\doteqdot 9298$ for $st(0.8)$. 

Consequently we notice that the fixed interval method could be extremely inefficient to the moving interval method if the unknown mother distribution assigns very small probability for one of the chosen intervals.  This could happen if the mother distribution has a finite support.  Suppose that we have prior knowledge that the mother distribution has the support $[0, 1]$, and set $a_i$'s as $a_i=i/10 (1\leq  i \leq 9)$ for the fixed intervals. The $m_i$'s for the moving interval method with "quantiles" are again $m_i=1/10 (0 \leq i \leq 9)$. 

\begin{figure}
\centering
\begin{subfigure}{0.4\columnwidth}
\centering
\includegraphics[width=\columnwidth]{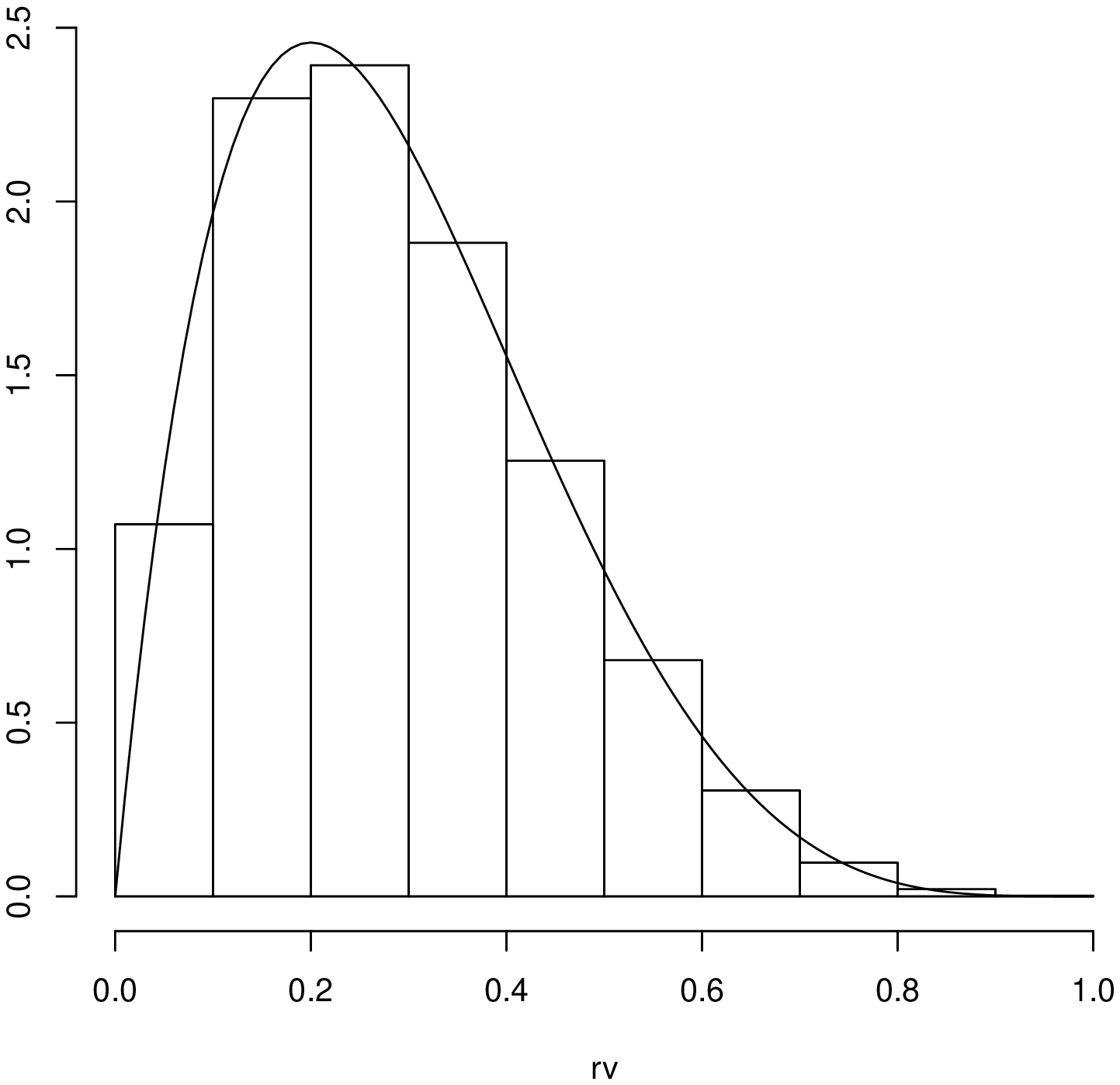}
\caption{Density and Histogram}
\label{dens_beta}
\end{subfigure}
\hspace{5mm}
\begin{subfigure}{0.4\columnwidth}
\centering
\includegraphics[width=\columnwidth]{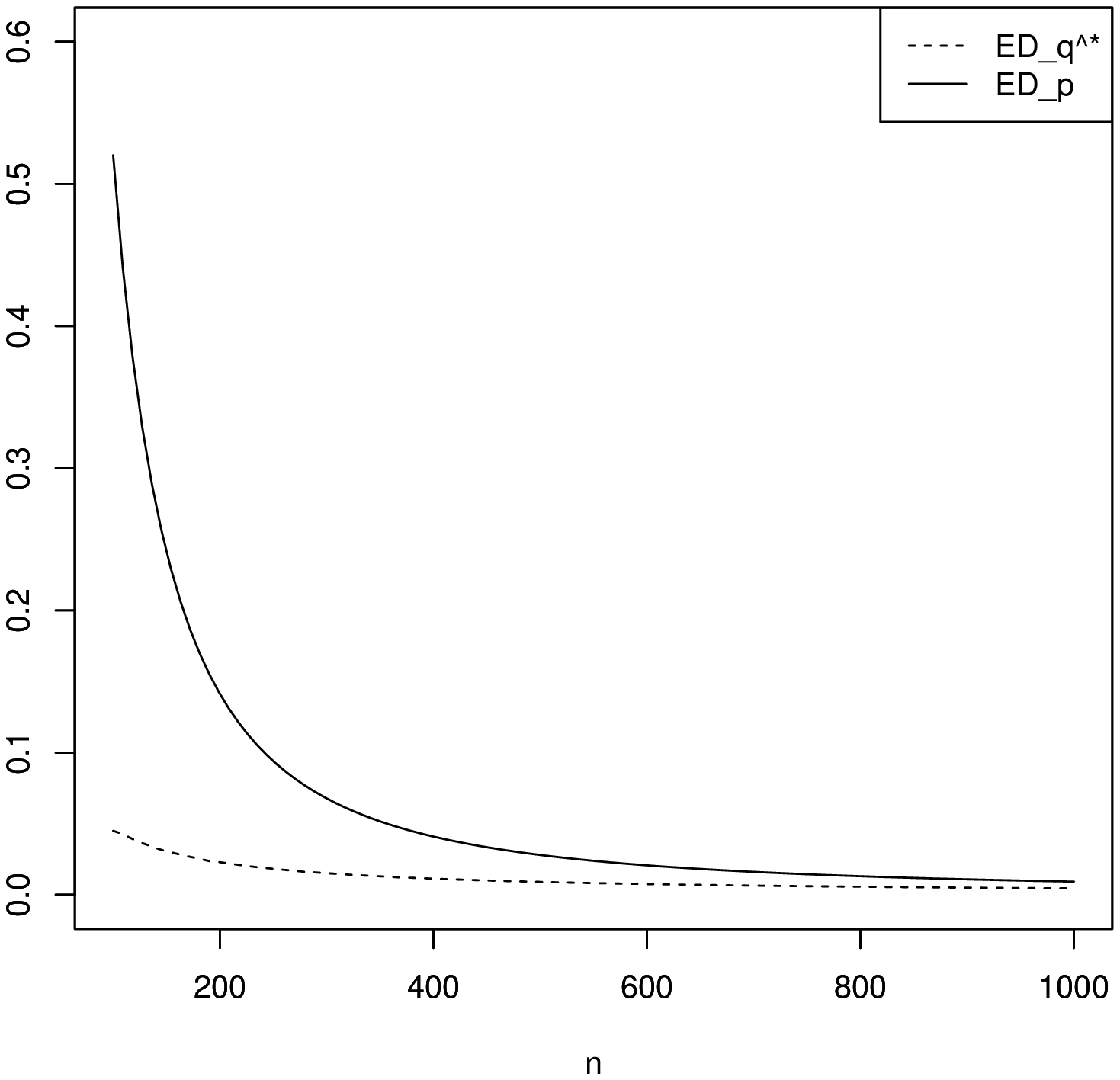}
\caption{The risks for the both methods}
\label{risks_beta}
\end{subfigure}
\caption{Beta(2,5)}
\label{risks_fin_suppo}
\end{figure}

If the mother distribution is $Beta(2, 5)$, the corresponding probabilities for the fixed intervals are given by
$$
(m_0, m_1, \ldots, m_9) \doteqdot (0.114, 0.230, 0.235, 0.187, 0.124, 0.068, 0.030, 0.009, 0.002, 5.5*10^{-5}).
$$
The graph of the density function and the histogram of $10^4$ samples with above $a_i$'s as the endpoints are given in Figure \ref{dens_beta} . The graphs of the approximated risks for the both methods are shown in Figure \ref{risks_beta}.
The solution for the equation \eqref{eq_EDI_EDP} is given by $n\doteqdot 379$. Even if we are lucky enough to know the finite support of the mother distribution, the fixed interval method is still quite inefficient to the moving interval method.

We saw that the moving interval method is superior theoretically and numerically to the fixed interval method as estimation of the mother distribution.  Needless to say, we often need to know the probability of some fixed intervals for a certain practical purpose. In that case, it might be preferable that the moving interval method is also subsidiarily used, since it could give some information on $M$ in \eqref{ED_I_expan} for the fixed interval method. Lastly we mention that the histogram (as estimation of the unknown distribution) falls between the both methods. In a conventional way, the intervals for the histogram are chosen after the sample is taken, taking into the consideration the frequency of each interval, especially being careful not to create the interval of null frequency.
\section{Appendix}
\label{append1}
--\textit{Proof of \eqref{E_R2}, \eqref{E_R3}, \eqref{E_R4}}--
\\
From \eqref{form_hat_m}, we notice that
$$
R_i=\frac{\hat{m}_i}{m_i}-1=\frac{1}{\sqrt{n}\:m_i}(\delip-\deli),
$$
hence
\begin{align}
\xwai m_i R_i^2 &=n^{-1} \xwai m_i^{-1} (\delip^2+\deli^2-2\deli\delip), \label{ex_form_E_R2}\\
\xwai m_i R_i^3 &=n^{-3/2} \xwai m_i^{-2}(\delip^3-3\delip^2\deli+3\delip\deli^2-\deli^3), \label{ex_form_E_R3}\\
\xwai m_i R_i^4 &=n^{-2} \xwai m_i^{-3}(\delip^4-4\delip^3\deli+6\delip^2\deli^2-4\delip\deli^3+\deli^4). \label{ex_form_E_R4}
\end{align}
From the formula on the moments of the ordered statistics $U_{(n_i)}$ (see (3.1.6) of \cite{David&Nagaraja})
\begin{equation}
E\Bigl[\prod_{i=1}^k U_{(n_i)}^{a_i} \Bigr]= \frac{n!}{\bigl(n+\sum_{i=1}^k a_i\bigr)!}\prod_{i=1}^k \frac{\bigl(n_i-1+\sum_{j=1}^i a_j \bigr)!}{\bigl(n_i-1+\sum_{j=1}^{i-1} a_j \bigr)!}, \qquad n_1 \leq \cdots \leq n_k,
\end{equation}
we have the following results.
\begin{align}
&E[U_{(\nui)}]\nonumber\\
&=\frac{n!}{(n+1)!}\frac{n_i !}{(n_i-1)!}\nonumber\\
&=\frac{n_i}{n+1} \nonumber\\
&=\frac{n}{n+1}\Bigl(\lami+\frac{\ri}{n}\Bigr) \nonumber\\
&=\Bigl(1-\frac{1}{n}+\frac{1}{n^2}+O(n^{-3})\Bigr)\Bigl(\lami+\frac{\ri}{n}\Bigr)\nonumber\\
&=\lami+\frac{1}{n}(-\lami+\ri)+\frac{1}{n^2}(-\ri+\lami)+O(n^{-3}),
\end{align}
where the forth equation is due to the fact
\begin{align*}
&\frac{n}{n+1}=1-\frac{1}{n+1}=1-\frac{1}{n}+\Bigl(\frac{1}{n}-\frac{1}{n+1}\Bigr)=1-\frac{1}{n}+\frac{1}{n(n+1)}\\
&1-\frac{1}{n}+\frac{1}{n^2}+\Bigl(-\frac{1}{n^2}+\frac{1}{n(n+1)}\Bigr)=1-\frac{1}{n}+\frac{1}{n^2}-\frac{1}{n^2(n+1)}=1-\frac{1}{n}+\frac{1}{n^2}+O(n^{-3}).
\end{align*}
\begin{align}
&E[U_{(\nui)}^2]\nonumber\\
&=\frac{n!}{(n+2)!}\frac{(n_i+1)!}{(n_i-1)!}\nonumber\\
&=\frac{1}{(n+1)(n+2)}n_i(n_i+1)\nonumber\\
&=\frac{n^2}{(n+1)(n+2)}\Bigl(\lami+\frac{\ri}{n}\Bigr)\Bigl(\lami+\frac{\ri+1}{n}\Bigr)\nonumber\\
&=\Bigl(1-\frac{3}{n}+\frac{7}{n^2}+O(n^{-3})\Bigr)\Bigl(\lami+\frac{\ri}{n}\Bigr)\Bigl(\lami+\frac{\ri+1}{n}\Bigr)\nonumber\\
&=\lami^2+\frac{1}{n}\bigl(-3\lami^2+\ri\lami+(\ri+1)\lami\bigr)\nonumber\\
&\quad+\frac{1}{n^2}\bigl(-3\ri\lami-3(\ri+1)\lami+\ri(\ri+1)+7\lami^2\bigr)+O(n^{-3})\nonumber\\
&=\lami^2+\frac{1}{n}\bigl(-3\lami^2+(2\ri+1)\lami\bigr)\nonumber\\
&\quad+\frac{1}{n^2}\bigl(-3(2\ri+1)\lami+\ri(\ri+1)+7\lami^2\bigr)+O(n^{-3}),
\end{align}
where the forth equation is due to the fact
\begin{align*}
&\frac{n^2}{(n+1)(n+2)}=\frac{n^2}{n^2+3n+2}=1-\frac{3n+2}{n^2+3n+2}=1-\frac{3}{n}+\Bigl(\frac{3}{n}-\frac{3n+2}{n^2+3n+2}\Bigr)\\
&=1-\frac{3}{n}+\frac{7n+6}{n^3+3n^2+2n}=1-\frac{3}{n}+\frac{7}{n^2}+\Bigl(-\frac{7}{n^2}+\frac{7n+6}{n^3+3n^2+2n}\Bigr)\\
&= 1-\frac{3}{n}+\frac{7}{n^2}+\frac{-7(n^2+3n+2)+7n^2+6n}{n^4+3n^3+2n^2}= 1-\frac{3}{n}+\frac{7}{n^2}+O(n^{-3}).
\end{align*}
\begin{align}
&E[U_{(\nui)}U_{(\nuip)}]\nonumber\\
&=\frac{n!}{(n+2)!}\frac{(\nui-1+1)!(\nuip-1+2)!}{(\nui-1)!(\nuip-1+1)!}\nonumber\\
&=\frac{\nui(\nuip+1)}{(n+1)(n+2)}\nonumber\\
&=\frac{n^2}{(n+1)(n+2)}\Bigl(\lami+\frac{\ri}{n}\Bigr)\Bigl(\lamip+\frac{\rip+1}{n}\Bigr)\nonumber\\
&=\Bigl(1-\frac{3}{n}+\frac{7}{n^2}+O(n^{-3})\Bigr)\Bigl(\lami+\frac{\ri}{n}\Bigr)\Bigl(\lamip+\frac{\rip+1}{n}\Bigr)\nonumber\\
&=\lami\lamip+\frac{1}{n}\Bigl(-3\lami\lamip+\ri\lamip+(\rip+1)\lami\Bigr)\nonumber\\
&\quad+\frac{1}{n^2}\bigl(-3\ri\lamip-3(\rip+1)\lami+\ri(\rip+1)+7\lami\lamip\bigr)+O(n^{-3}).\nonumber\\
\end{align}
\begin{align}
&E[U_{(\nui)}^3]\nonumber\\
&=\frac{n!}{(n+3)!}\frac{(\nui-1+3)!}{(\nui-1)!}\nonumber\\
&=\frac{\nui(\nui+1)(\nui+2)}{(n+1)(n+2)(n+3)}\nonumber\\
&=\frac{n^3}{(n+1)(n+2)(n+3)}\Bigl(\lami+\frac{\ri}{n}\Bigr)\Bigl(\lami+\frac{\ri+1}{n}\Bigr)\Bigl(\lami+\frac{\ri+2}{n}\Bigr)\nonumber\\
&=\Bigl(1-\frac{6}{n}+\frac{25}{n^2}+O(n^{-3})\Bigr)\Bigl(\lami+\frac{\ri}{n}\Bigr)\Bigl(\lami+\frac{\ri+1}{n}\Bigr)\Bigl(\lami+\frac{\ri+2}{n}\Bigr)\nonumber\\
&=\lami^3+\frac{1}{n}\Bigl(-6\lami^3+\lami^2 \ri+\lami^2(\ri+1)+\lami^2(\ri+2)\Bigr)\nonumber\\
&\quad+\frac{1}{n^2}\Bigl(-6\ri\lami^2-6(\ri+1)\lami^2-6(\ri+2)\lami^2+\ri(\ri+1)\lami+\ri(\ri+2)\lami \nonumber\\
&\qquad\qquad+(\ri+1)(\ri+2)\lami+25\lami^3\Bigr)+O(n^{-3})\nonumber\\
&=\lami^3+\frac{1}{n}\Bigl(-6\lami^3+3\ri\lami^2+3\lami^2\Bigr)\nonumber\\
&\quad+\frac{1}{n^2}\Bigl(25\lami^3+(-18\ri-18)\lami^2+(3\ri^2+6\ri+2)\lami\Bigr)+O(n^{-3}),
\end{align}
where the forth equation is due to the following relation;
\begin{align}
&\frac{n^3}{(n+1)(n+2)(n+3)}-1=\frac{n^3-(n^2+3n+2)(n+3)}{(n+1)(n+2)(n+3)}\nonumber\\
&=\frac{-6n^2-11n-6}{(n+1)(n+2)(n+3)}=-\frac{6}{n}+\frac{6}{n}-\frac{6n^2+11n+6}{(n+1)(n+2)(n+3)}\nonumber\\
&=-\frac{6}{n}+\frac{25n^2+60n+36}{n(n^3+6n^2+11n+6)}\nonumber\\
&=-\frac{6}{n}+\frac{25}{n^2}+\Bigl(-\frac{25}{n^2}+\frac{25n^2+60n+36}{n(n^3+6n^2+11n+6)}\Bigr)\nonumber\\
&=-\frac{6}{n}+\frac{25}{n^2}+\frac{-25(n^3+6n^2+11n+6)+25n^3+60n^2+36n}{n^2(n^3+6n^2+11n+6)}\nonumber\\
&=-\frac{6}{n}+\frac{25}{n^2}+O(n^{-3}).\nonumber
\end{align}
\begin{align}
&E[U_{(\nui)}^2U_{(\nuip)}]\nonumber\\
&=\frac{n!}{(n+3)!}\frac{(\nui+1)!\:(\nuip+2)!}{(\nui-1)!\:(\nuip+1)!}\nonumber\\
&=\frac{\nui(\nui+1)(\nuip+2)}{(n+1)(n+2)(n+3)}\nonumber\\
&=\frac{n^3}{(n+1)(n+2)(n+3)}\Bigl(\lami+\frac{\ri}{n}\Bigr)\Bigl(\lami+\frac{\ri+1}{n}\Bigr)\Bigl(\lamip+\frac{\rip+2}{n}\Bigr)\nonumber\\
&=\Bigl(1-\frac{6}{n}+\frac{25}{n^2}\Bigr)\Bigl(\lami+\frac{\ri}{n}\Bigr)\Bigl(\lami+\frac{\ri+1}{n}\Bigr)\Bigl(\lamip+\frac{\rip+2}{n}\Bigr)\nonumber\\
&=\lami^2\lamip+\frac{1}{n}\Bigl(-6\lami^2\lamip+\ri\lami\lamip+(\ri+1)\lami\lamip+(\rip+2)\lami^2\Bigr)\nonumber\\
&\quad+\frac{1}{n^2}\Bigl(-6\ri\lami\lamip-6(\ri+1)\lami\lamip-6(\rip+2)\lami^2+\ri(\ri+1)\lamip\nonumber\\
&\qquad\qquad+\ri(\rip+2)\lami+(\ri+1)(\rip+2)\lami+25\lami^2\lamip\Bigr)+O(n^{-3})\nonumber\\
&=\lami^2\lamip+\frac{1}{n}\Bigl(-6\lami^2\lamip+(2\ri+1)\lami\lamip+(\rip+2)\lami^2\Bigr)\nonumber\\
&\quad+\frac{1}{n^2}\Bigl(25\lami^2\lamip-(6\rip+12)\lami^2-(12\ri+6)\lami\lamip\nonumber\\
&\qquad\qquad+(2\ri\rip+4\ri+\rip+2)\lami+(\ri^2+\ri)\lamip\Bigr)+O(n^{-3}).
\end{align}
\begin{align}
&E[U_{(\nui)}U_{(\nuip)}^2]\nonumber\\
&=\frac{n!}{(n+3)!}\frac{\nui!\:(\nuip+2)!}{(\nui-1)!\:\nuip!}\nonumber\\
&=\Bigl(1-\frac{6}{n}+\frac{25}{n^2}+O(n^{-3})\Bigr)\Bigl(\lami+\frac{\ri}{n}\Bigr)\Bigl(\lamip+\frac{\rip+1}{n}\Bigr)\Bigl(\lamip+\frac{\rip+2}{n}\Bigr)\nonumber\\
&=\lamip^2\lami+\frac{1}{n}\Bigl(-6\lamip^2\lami+\ri\lamip^2+(\rip+1)\lami\lamip+(\rip+2)\lami\lamip\Bigr)\nonumber\\
&\quad+\frac{1}{n^2}\Bigl(-6\ri\lamip^2-6(\rip+1)\lami\lamip-6(\rip+2)\lami\lamip+\ri(\rip+1)\lamip\nonumber\\
&\qquad\qquad+\ri(\rip+2)\lamip+(\rip+1)(\rip+2)\lami+25\lamip^2\lami\Bigr)+O(n^{-3})\nonumber\\
&=\lamip^2\lami+\frac{1}{n}\Bigl(-6\lamip^2\lami+\ri\lamip^2+(2\rip+3)\lami\lamip\Bigr)\nonumber\\
&\quad+\frac{1}{n^2}\Bigl(25\lamip^2\lami-6\ri\lamip^2-(12\rip+18)\lami\lamip\nonumber\\
&\qquad\qquad+(\rip^2+3\rip+2)\lami+(2\ri\rip+3\ri)\lamip\Bigr)+O(n^{-3}).
\end{align}
\begin{align}
&E[U_{(\nui)}^4]\nonumber\\
&=\frac{n!}{(n+4)!}\frac{(\nui+3)!}{(\nui-1)!}\nonumber\\
&=\frac{n^4}{(n+1)(n+2)(n+3)(n+4)}\frac{\nui(\nui+1)(\nui+2)(\nui+3)}{n^4}\nonumber\\
&=\Bigl(1-\frac{10}{n}+\frac{65}{n^2}+O(n^{-3})\Bigr)\Bigl(\lami+\frac{\ri}{n}\Bigr)\Bigl(\lami+\frac{\ri+1}{n}\Bigr)\Bigl(\lami+\frac{\ri+2}{n}\Bigr)\Bigl(\lami+\frac{\ri+3}{n}\Bigr)\nonumber\\
&=\lami^4+\frac{1}{n}\Bigl(-10\lami^4+\ri\lami^3+(\ri+1)\lami^3+(\ri+2)\lami^3+(\ri+3)\lami^3\Bigr)\nonumber\\
&\quad+\frac{1}{n^2}\Bigl(-10\ri\lami^3-10(\ri+1)\lami^3-10(\ri+2)\lami^3-10(\ri+3)\lami^3+\ri(\ri+1)\lami^2\nonumber\\
&\qquad\qquad+\ri(\ri+2)\lami^2+\ri(\ri+3)\lami^2+(\ri+1)(\ri+2)\lami^2\nonumber\\
&\qquad\qquad+(\ri+1)(\ri+3)\lami^2+(\ri+2)(\ri+3)\lami^2+65\lami^4\Bigr)+O(n^{-3})\nonumber\\
&=\lami^4+\frac{1}{n}\Bigl(-10\lami^4+(4\ri+6)\lami^3\Bigr)\nonumber\\
&\quad+\frac{1}{n^2}\Bigl( 65\lami^4-(40\ri+60)\lami^3+(6\ri^2+18\ri+11)\lami^2\Bigr)+O(n^{-3}),
\end{align}
where the third equation is due to the following relation;
\begin{align}
&\frac{n^4}{(n+1)(n+2)(n+3)(n+4)}-1=\frac{n^4-(n^4+10n^3+35n^2+50n+24)}{(n+1)(n+2)(n+3)(n+4)}\nonumber\\
&=\frac{-(10n^3+35n^2+50n+24)}{n^4+10n^3+35n^2+50n+24}=-\frac{10}{n}+\frac{10}{n}-\frac{10n^3+35n^2+50n+24}{n^4+10n^3+35n^2+50n+24}\nonumber\\
&=-\frac{10}{n}+\frac{10n^4+100n^3+350n^2+500n+240-10n^4-35n^3-50n^2-24n}{n^5+10n^4+35n^3+50n^2+24n}\nonumber\\
&=-\frac{10}{n}+\frac{65n^3+300n^2+476n+240}{n^5+10n^4+35n^3+50n^2+24n}\nonumber\\
&=-\frac{10}{n}+\frac{65}{n^2}+O(n^{-3}).
\end{align}
\begin{align}
&E[U_{(\nui)}^3U_{(\nuip)}]\nonumber\\
&=\frac{n!}{(n+4)!}\frac{(\nui+2)!\:(\nuip+3)!}{(\nui-1)!\:(\nuip+2)!}\nonumber\\
&=\frac{n^4}{(n+1)(n+2)(n+3)(n+4)}\frac{\nui(\nui+1)(\nui+2)(\nuip+3)}{n^4}\nonumber\\
&=\Bigl(1-\frac{10}{n}+\frac{65}{n^2}+O(n^{-3})\Bigr)\Bigl(\lami+\frac{\ri}{n}\Bigr)\Bigl(\lami+\frac{\ri+1}{n}\Bigr)\Bigl(\lami+\frac{\ri+2}{n}\Bigr)\nonumber\\
&\qquad\times\Bigl(\lamip+\frac{\rip+3}{n}\Bigr)\nonumber\\
&=\lami^3\lamip+\frac{1}{n}\Bigl(-10\lami^3\lamip+\ri\lami^2\lamip+(\ri+1)\lami^2\lamip+(\ri+2)\lami^2\lamip\nonumber\\
&\hspace{30mm}+(\rip+3)\lami^3\Bigr)\nonumber\\
&\quad+\frac{1}{n^2}\Bigl(-10\ri\lami^2\lamip-10(\ri+1)\lami^2\lamip-10(\ri+2)\lami^2\lamip-10(\rip+3)\lami^3\nonumber\\
&\qquad\qquad+\ri(\ri+1)\lami\lamip+\ri(\ri+2)\lami\lamip+\ri(\rip+3)\lami^2\nonumber\\
&\qquad\qquad+(\ri+1)(\ri+2)\lami\lamip+(\ri+1)(\rip+3)\lami^2+(\ri+2)(\rip+3)\lami^2\nonumber\\
&\qquad\qquad+65\lami^3\lamip\Bigr)+O(n^{-3})\nonumber\\
&=\lami^3\lamip+\frac{1}{n}\Bigl(-10\lami^3\lamip+(3\ri+3)\lami^2\lamip+(\rip+3)\lami^3\Bigr)\nonumber\\
&\quad+\frac{1}{n^2}\Bigl(65\lami^3\lamip-10(\rip+3)\lami^3-(30\ri+30)\lami^2\lamip\nonumber\\
&\qquad \qquad+(3\ri\rip+9\ri+3\rip+9)\lami^2+(3\ri^2+6\ri+2)\lami\lamip\Bigr)+O(n^{-3}).
\end{align}
\begin{align}
&E[U_{(\nui)}U_{(\nuip)}^3]\nonumber\\
&=\frac{n!}{(n+4)!}\frac{\nui!\:(\nuip+3)!}{(\nui-1)!\:\nuip!}\nonumber\\
&=\frac{n^4}{(n+1)(n+2)(n+3)(n+4)}\frac{\nui(\nuip+1)(\nuip+2)(\nuip+3)}{n^4}\nonumber\\
&=\Bigl(1-\frac{10}{n}+\frac{65}{n^2}+O(n^{-3})\Bigr)\Bigl(\lami+\frac{\ri}{n}\Bigr)\Bigl(\lamip+\frac{\rip+1}{n}\Bigr)\Bigl(\lamip+\frac{\rip+2}{n}\Bigr)\nonumber\\
&\qquad\times\Bigl(\lamip+\frac{\rip+3}{n}\Bigr)\nonumber\\
&=\lami\lamip^3+\frac{1}{n}\Bigl(-10\lami\lamip^3+\ri\lamip^3+(\rip+1)\lami\lamip^2+(\rip+2)\lami\lamip^2\nonumber\\
&\hspace{30mm}+(\rip+3)\lami\lamip^2\Bigr)\nonumber\\
&\quad+\frac{1}{n^2}\Bigl(-10\ri\lamip^3-10(\rip+1)\lami\lamip^2-10(\rip+2)\lami\lamip^2-10(\rip+3)\lami\lamip^2\nonumber\\
&\qquad\qquad+\ri(\rip+1)\lamip^2+\ri(\rip+2)\lamip^2+\ri(\rip+3)\lamip^2\nonumber\\
&\qquad\qquad+(\rip+1)(\rip+2)\lami\lamip+(\rip+1)(\rip+3)\lami\lamip\nonumber\\
&\qquad\qquad+(\rip+2)(\rip+3)\lami\lamip+65\lami\lamip^3\Bigr)+O(n^{-3})\nonumber\\
&=\lami\lamip^3+\frac{1}{n}\Bigl(-10\lami\lamip^3+\ri\lamip^3+(3\rip+6)\lami\lamip^2\Bigr)\nonumber\\
&\quad+\frac{1}{n^2}\Bigl(65\lami\lamip^3-10\ri\lamip^3-(30\rip+60)\lami\lamip^2+(3\ri\rip+6\ri)\lamip^2\nonumber\\
&\qquad\qquad+(3\rip^2+12\rip+11)\lami\lamip\Bigr)+O(n^{-3}).
\end{align}
\begin{align}
&E[U_{(\nui)}^2U_{(\nuip)}^2]\nonumber\\
&=\frac{n!}{(n+4)!}\frac{(\nui+1)!\:(\nuip+3)!}{(\nui-1)!\:(\nuip+1)!}\nonumber\\
&=\frac{n^4}{(n+1)(n+2)(n+3)(n+4)}\frac{\nui(\nui+1)(\nuip+2)(\nuip+3)}{n^4}\nonumber\\
&=\Bigl(1-\frac{10}{n}+\frac{65}{n^2}+O(n^{-3})\Bigr)\Bigl(\lami+\frac{\ri}{n}\Bigr)\Bigl(\lami+\frac{\ri+1}{n}\Bigr)\Bigl(\lamip+\frac{\rip+2}{n}\Bigr)\nonumber\\
&\qquad\times\Bigl(\lamip+\frac{\rip+3}{n}\Bigr)\nonumber\\
&=\lami^2\lamip^2+\frac{1}{n}\Bigl(-10\lami^2\lamip^2+\ri\lami\lamip^2+(\ri+1)\lami\lamip^2+(\rip+2)\lami^2\lamip\nonumber\\
&\hspace{30mm}+(\rip+3)\lami^2\lamip\Bigr)\nonumber\\
&\quad+\frac{1}{n^2}\Bigl(-10\ri\lami\lamip^2-10(\ri+1)\lami\lamip^2-10(\rip+2)\lami^2\lamip-10(\rip+3)\lami^2\lamip\nonumber\\
&\qquad\qquad+\ri(\ri+1)\lamip^2+\ri(\rip+2)\lami\lamip+\ri(\rip+3)\lami\lamip\nonumber\\
&\qquad\qquad+(\ri+1)(\rip+2)\lami\lamip+(\ri+1)(\rip+3)\lami\lamip\nonumber\\
&\qquad\qquad+(\rip+2)(\rip+3)\lami^2+65\lami^2\lamip^2\Bigr)+O(n^{-3})\nonumber\\
&=\lami^2\lamip^2+\frac{1}{n}\Bigl(-10\lami^2\lamip^2+(2\rip+5)\lami^2\lamip+(2\ri+1)\lami\lamip^2\Bigr)\nonumber\\
&\quad+\frac{1}{n^2}\Bigl(65\lami^2\lamip^2-(20\ri+10)\lami\lamip^2-(20\rip+50)\lami^2\lamip+(\rip^2+5\rip+6)\lami^2\nonumber\\
&\qquad\qquad+(\ri^2+\ri)\lamip^2+(4\ri\rip+10\ri+2\rip+5)\lami\lamip\Bigr)+O(n^{-3}).
\end{align}
From the moments of $\uni$'s, we can calculate the moments of $\deli$'s as follows.
\begin{align}
&n^{-1}E[\deli^2]\nonumber\\
&=E[(\uni-\lami)^2]\nonumber\\
&=E[\uni^2-2\lami \uni+\lami^2]\nonumber\\
&=\lami^2-2\lami^2+\lami^2+n^{-1}\bigl(-3\lami^2+(2\ri+1)\lami+2\lami^2-2\ri\lami\bigr)\nonumber\\
&\qquad+n^{-2}\bigl(-3(2\ri+1)\lami+\ri(\ri+1)+7\lami^2+2\ri\lami-2\lami^2\bigr)+O(n^{-3})\nonumber\\
&=\frac{1}{n}\lami(1-\lami)+\frac{1}{n^2}\bigl(5\lami^2-3\lami-4\ri\lami+\ri(\ri+1)\bigr)+O(n^{-3}).
\end{align}
\begin{align}
&n^{-1}E[\deli\delip]\nonumber\\
&=E[(\uni-\lami)(\unip-\lamip)]\nonumber\\
&=E[\uni\unip-\lamip \uni-\lami \unip+\lami\lamip]\nonumber\\
&=\lami\lamip+n^{-1}\bigl(-3\lami\lamip+\ri\lamip+\rip\lami+\lami\bigr)\nonumber\\
&\quad+n^{-2}\bigl(-3\ri\lamip-3(\rip+1)\lami+7\lami\lamip+\ri(\rip+1)\bigr)\nonumber\\
&\quad-\lami\lamip-n^{-1}\bigl(-\lami\lamip+\lami\rip\bigr)-n^{-2}\bigl(-\rip\lami+\lamip\lami\bigr)\nonumber\\
&\quad-\lami\lamip-n^{-1}\bigl(-\lami\lamip+\lamip\ri\bigr)-n^{-2}\bigl(-\ri\lamip+\lamip\lami\bigr)\nonumber\\
&\quad+\lami\lamip+O(n^{-3})\nonumber\\
&=n^{-1}\bigl(-\lami\lamip+\lami\bigr)+n^{-2}\bigl(5\lami\lamip-2\ri\lamip-2\rip\lami-3\lami+\ri(\rip+1)\bigr)+O(n^{-3}).
\end{align}
\begin{align}
&n^{-3/2}E[\deli^3]\nonumber\\
&=E[(\uni-\lami)^3]\nonumber\\
&=E[\uni^3-3\lami\uni^2+3\lami^2\uni-\lami^3]\nonumber\\
&=\lami^3+n^{-1}\bigl(-6\lami^3+3\ri\lami^2+3\lami^2\bigr)\nonumber\\
&\quad+n^{-2}\bigl(25\lami^3-(18\ri+18)\lami^2+(3\ri^2+6\ri+2)\lami\bigr)\nonumber\\
&\quad-3\lami^3+n^{-1}\bigl(9\lami^3-6\ri\lami^2-3\lami^2\bigr)\nonumber\\
&\quad+n^{-2}\bigl(18\ri\lami^2+9\lami^2-3\ri(\ri+1)\lami-21\lami^3\bigr)\nonumber\\
&\quad+3\lami^3+n^{-1}\bigl(-3\lami^3+3\ri\lami^2\bigr)+n^{-2}\bigl(-3\ri\lami^2+3\lami^3\bigr)-\lami^3
+O(n^{-3})\nonumber\\
&=n^{-2}\bigl(7\lami^3-3\lami^2(\ri+3)+\lami(3\ri+2)\bigr)+O(n^{-3}).
\end{align}
\begin{align}
&n^{-3/2}E[\deli^2\delip]\nonumber\\
&=E[(\uni-\lami)^2(\unip-\lamip)]\nonumber\\
&=E[\uni^2\unip-\lamip\uni^2-2\lami\uni\unip+2\lami\lamip\uni+\lami^2\unip-\lami^2\lamip]\nonumber\\
&=\lami^2\lamip+n^{-1}\bigl(-6\lami^2\lamip+(2\ri+1)\lami\lamip+(\rip+2)\lami^2\bigr)\nonumber\\
&\quad+n^{-2}\bigl(25\lami^2\lamip-(6\rip+12)\lami^2-(12\ri+6)\lami\lamip\nonumber\\
&\qquad\qquad +(2\ri\rip+4\ri+\rip+2)\lami+(\ri^2+\ri)\lamip\bigr)\nonumber\\
&\quad-\lami^2\lamip+n^{-1}\bigl(3\lami^2\lamip-(2\ri+1)\lami\lamip\bigr)\nonumber\\
&\quad+n^{-2}\bigl(3(2\ri+1)\lami\lamip-\ri(\ri+1)\lamip-7\lami^2\lamip\bigr)\nonumber\\
&\quad-2\lami^2\lamip+n^{-1}\bigl(6\lami^2\lamip-2\ri\lami\lamip-2(\rip+1)\lami^2\bigr)\nonumber\\
&\quad+n^{-2}\bigl(6\ri\lami\lamip+6(\rip+1)\lami^2-2\ri(\rip+1)\lami-14\lami^2\lamip\bigr)\nonumber\\
&\quad+2\lami^2\lamip+n^{-1}\bigl(-2\lami^2\lamip+2\lami\lamip\ri\bigr)+n^{-2}\bigl(-2\lami\lamip\ri+2\lami^2\lamip\bigr)\nonumber\\
&\quad+\lami^2\lamip+n^{-1}\bigl(-\lami^2\lamip+\rip\lami^2\bigr)+n^{-2}\bigl(-\lami^2\rip+\lami^2\lamip\bigr)-\lami^2\lamip+O(n^{-3})\nonumber\\
&=n^{-2}\bigl(7\lami^2\lamip-(\rip+6)\lami^2-(2\ri+3)\lami\lamip+(2\ri+\rip+2)\lami\bigr)+O(n^{-3}).
\end{align}
\begin{align}
&n^{-3/2}E[\deli\delip^2]\nonumber\\
&=E[(\uni-\lami)(\unip-\lamip)^2]\nonumber\\
&=E[\unip^2\uni-\lami\unip^2-2\lamip\uni\unip+2\lami\lamip\unip+\lamip^2\uni-\lamip^2\lami]\nonumber\\
&=\lamip^2\lami+n^{-1}\bigl(-6\lamip^2\lami+\ri\lamip^2+(2\rip+3)\lami\lamip\bigr)\nonumber\\\
&\quad+n^{-2}\bigl(-6\ri\lamip^2-(12\rip+18)\lami\lamip+(\rip^2+3\rip+2)\lami\nonumber\\
&\quad\qquad\quad+(2\ri\rip+3\rip)\lamip+25\lamip^2\lami\bigr)\nonumber\\
&\quad-\lamip^2\lami+n^{-1}\bigl(3\lamip^2\lami-(2\rip+1)\lami\lamip \bigr)\nonumber\\
&\quad+n^{-2}\bigl(6\rip\lamip\lami+3\lami\lamip-\rip(\rip+1)\lami-7\lamip^2\lami\bigr)\nonumber\\
&\quad-2\lami\lamip^2+n^{-1}\bigl(6\lami\lamip^2-2\ri\lamip^2-2(\rip+1)\lami\lamip\bigr)\nonumber\\
&\quad+n^{-2}\bigl(6\ri\lamip^2+6(\rip+1)\lami\lamip-2\ri(\rip+1)\lamip-14\lami\lamip^2\bigr)\nonumber\\
&\quad+2\lami\lamip^2+n^{-1}\bigl(-2\lami\lamip^2+2\rip\lami\lamip\bigr)+n^{-2}\bigl(-2\rip\lami\lamip+2\lami\lamip^2\bigr)\nonumber\\
&\quad+\lami\lamip^2+n^{-1}\bigl(-\lami\lamip^2+\ri\lamip^2\bigr)+n^{-2}\bigl(-\ri\lamip^2+\lami\lamip^2\bigr)-\lami\lamip^2+O(n^{-3})\nonumber\\
&=n^{-2}\bigl(\lamip^2\lami(25-7-14+2+1)+\lamip^2(-6\ri+6\ri-\ri)\nonumber\\
&\qquad\quad+\lami\lamip(-12\rip-18+6\rip+3+6\rip+6-2\rip)\nonumber\\
&\qquad\quad+\lamip(2\ri\rip+3\ri-2\ri\rip-2\ri)\nonumber\\
&\qquad\quad+\lami(\rip^2+3\rip+2-\rip^2-\rip)\bigr)\nonumber\\
&=n^{-2}\bigl(7\lamip^2\lami-\ri\lamip^2-(2\rip+9)\lami\lamip+\ri\lamip+(2\rip+2)\lami\bigr)+O(n^{-3}).
\end{align}
\begin{align}
&n^{-2}E[\deli^4]\nonumber\\
&=E(\uni-\lami)^4]\nonumber\\
&=E[\uni^4-4\lami\uni^3+6\lami^2\uni^2-4\lami^3\uni+\lami^4]\nonumber\\
&=\lami^4+n^{-1}\bigl(-10\lami^4+(4\ri+6)\lami^3\bigr)+n^{-2}\bigl(65\lami^4+(-40\ri-60)\lami^3+(6\ri^2+18\ri+11)\lami^2\bigr)\nonumber\\
&\quad-4\lami^4+n^{-1}\bigl(24\lami^4-12\lami^3\ri-12\lami^3\bigr)\nonumber\\
&\quad+n^{-2}\bigl(-100\lami^4+(72\ri+72)\lami^3+(-12\ri^2-24\ri-8)\lami^2\bigr)\nonumber\\
&\quad+6\lami^4+n^{-1}\bigl(-18\lami^4+(12\ri+6)\lami^3\bigr)+n^{-2}\bigl(42\lami^4+(-36\ri-18)\lami^3+6(\ri^2+\ri)\lami^2\bigr)\nonumber\\
&\quad-4\lami^4+n^{-1}\bigl(4\lami^4-4\lami^3\ri\bigr)+n^{-2}\bigl(4\ri\lami^3-4\lami^4\bigr)+\lami^4+O(n^{-3})\nonumber\\
&=n^{-2}\bigl(3\lami^4-6\lami^3+3\lami^2\bigr)+O(n^{-3}).
\end{align}
\begin{align}
&n^{-2}E[\deli^3\delip]\nonumber\\
&=E[(\uni-\lami)^3(\unip-\lamip)]\nonumber\\
&=E[(\uni^3-3\uni^2\lami+3\uni\lami^2-\lami^3)(\unip-\lamip)]\nonumber\\
&=E[\uni^3\unip-\lamip\uni^3-3\lami\uni^2\unip+3\lami\lamip\uni^2+3\lami^2\uni\unip\nonumber\\
&\qquad-3\lami^2\lamip\uni-\lami^3\unip+\lami^3\lamip]\nonumber\\
&=\lami^3\lamip+n^{-1}\bigl( -10\lami^3\lamip+(\rip+3)\lami^3+(3\ri+3)\lami^2\lamip\bigr)\nonumber\\
&\quad+n^{-2}\bigl(65\lami^3\lamip+(-10\rip-30)\lami^3+(-30\ri-30)\lami^2\lamip\nonumber\\
&\qquad\qquad+(3\ri\rip+9\ri+3\rip+9)\lami^2+(3\ri^2+6\ri+2)\lami\lamip\bigr)\nonumber\\
&\quad-\lami^3\lamip+n^{-1}\bigl(6\lami^3\lamip+(-3\ri-3)\lami^2\lamip\bigr)\nonumber\\
&\quad+n^{-2}\bigl(-25\lami^3\lamip+(18\ri+18)\lami^2\lamip+(-3\ri^2-6\ri-2)\lami\lamip\bigr)\nonumber\\
&\quad-3\lami^3\lamip+n^{-1}\bigl(18\lami^3\lamip+(-3\rip-6)\lami^3+(-6\ri-3)\lami^2\lamip\bigr)\nonumber\\
&\quad+n^{-2}\bigl(-75\lami^3\lamip+(18\rip+36)\lami^3+(-6\ri\rip-12\ri-3\rip-6)\lami^2\nonumber\\
&\qquad\qquad+(36\ri+18)\lami^2\lamip+(-3\ri^2-3\ri)\lami\lamip\bigr)\nonumber\\
&\quad+3\lami^3\lamip+n^{-1}\bigl(-9\lami^3\lamip+(6\ri+3)\lami^2\lamip\bigr)\nonumber\\
&\quad+n^{-2}\bigl(21\lami^3\lamip+(-18\ri-9)\lami^2\lamip+(3\ri^2+3\ri)\lami\lamip\bigr)\nonumber\\
&\quad+3\lami^3\lamip+n^{-1}\bigl(-9\lami^3\lamip+(3\rip+3)\lami^3+3\ri\lami^2\lamip\bigr)\nonumber\\
&\quad+n^{-2}\bigl(21\lami^3\lamip+(-9\rip-9)\lami^3+(3\ri\rip+3\ri)\lami^2-9\ri\lami^2\lamip\bigr)\nonumber\\
&\quad-3\lami^3\lamip+n^{-1}\bigl(3\lami^3\lamip-3\ri\lami^2\lamip\bigr)+n^{-2}\bigl(-3\lami^3\lamip+3\ri\lami^2\lamip\bigr)\nonumber\\
&\quad-\lami^3\lamip+n^{-1}\bigl(\lami^3\lamip-\rip\lami^3\bigr)+n^{-2}\bigl(-\lami^3\lamip+\rip\lami^3\bigr)+\lami^3\lamip+O(n^{-3})\nonumber\\
&=n^{-2}\bigl(3\lami^3\lamip-3\lami^3-3\lami^2\lamip+3\lami^2\bigr)+O(n^{-3}).
\end{align}
\begin{align}
&n^{-2}E[\deli\delip^3]\nonumber\\
&=E[(\unip-\lamip)^3(\uni-\lami)]\nonumber\\
&=E[(\unip^3-3\unip^2\lamip+3\unip\lamip^2-\lamip^3)(\uni-\lami)]\nonumber\\
&=E[\unip^3\uni-\lami\unip^3-3\lamip\unip^2\uni+3\lamip\lami\unip^2+3\lamip^2\unip\uni\nonumber\\
&\qquad-3\lamip^2\lami\unip-\lamip^3\uni+\lamip^3\lami]\nonumber\\
&=\lami\lamip^3+n^{-1}\bigl( -10\lami\lamip^3+\ri\lamip^3+(3\rip+6)\lami\lamip^2\bigr)\nonumber\\
&\quad+n^{-2}\bigl(65\lami\lamip^3-10\ri\lamip^3-(30\rip+60)\lami\lamip^2\nonumber\\
&\qquad\qquad+(3\ri\rip+6\ri)\lamip^2+(3\rip^2+12\rip+11)\lami\lamip\bigr)\nonumber\\
&\quad-\lami\lamip^3+n^{-1}\bigl(6\lami\lamip^3-(3\rip+3)\lamip^2\lami\bigr)\nonumber\\
&\quad+n^{-2}\bigl((18\rip+18)\lami\lamip^2-(3\rip^2+6\rip+2)\lami\lamip-25\lamip^3\lami\bigr)\nonumber\\
&\quad-3\lami\lamip^3+n^{-1}\bigl(18\lamip^3\lami-3\ri\lamip^3-3(2\rip+3)\lamip^2\lami\bigr)\nonumber\\
&\quad+n^{-2}\bigl(-75\lami\lamip^3+18\ri\lamip^3+(36\rip+54)\lami\lamip^2\nonumber\\
&\qquad\qquad-(3\rip^2+9\rip+6)\lami\lamip-(6\ri\rip+9\ri)\lamip^2\bigr)\nonumber\\
&\quad+3\lami\lamip^3+n^{-1}\bigl(-9\lami\lamip^3+(6\rip+3)\lami\lamip^2\bigr)\nonumber\\
&\quad+n^{-2}\bigl(21\lamip^3\lami-(18\rip+9)\lamip^2\lami+(3\rip^2+3\rip)\lami\lamip\bigr)\nonumber\\
&\quad+3\lami\lamip^3+n^{-1}\bigl(-9\lami\lamip^3+3\ri\lamip^3+(3\rip+3)\lami\lamip^2\bigr)\nonumber\\
&\quad+n^{-2}\bigl(21\lami\lamip^3-9\ri\lamip^3-(9\rip+9)\lami\lamip^2+(3\ri\rip+3\ri)\lamip^2\bigr)\nonumber\\
&\quad-3\lami\lamip^3+n^{-1}\bigl(3\lamip^3\lami-3\rip\lamip^2\lami\bigr)+n^{-2}\bigl(3\rip\lamip^2\lami-3\lamip^3\lami\bigr)\nonumber\\
&\quad-\lami\lamip^3+n^{-1}\bigl(\lamip^3\lami-\ri\lamip^3\bigr)+n^{-2}\bigl(\ri\lamip^3-\lamip^3\lami\bigr)+\lamip^3\lami+O(n^{-3})\nonumber\\
&=n^{-2}\bigl(3\lami\lamip^3-6\lamip^2\lami+3\lami\lamip\bigr)+O(n^{-3}).
\end{align}
\begin{align}
&n^{-2}E[\deli^2\delip^2]\nonumber\\
&=E[(\unip-\lamip)^2(\uni-\lami)^2]\nonumber\\
&=E[(\uni^2-2\lami\uni+\lami^2)(\unip^2-2\lamip\unip+\lamip^2)]\nonumber\\
&=E[\uni^2\unip^2-2\lamip\uni^2\unip+\lamip^2\uni^2-2\lami\uni\unip^2\nonumber\\
&\qquad+4\lami\lamip\uni\unip-2\lami\lamip^2\uni+\lami^2\unip^2-2\lami^2\lamip\unip+\lami^2\lamip^2]\nonumber\\
&=\lami^2\lamip^2+n^{-1}\bigl(-10\lami^2\lamip^2+(2\rip+5)\lami^2\lamip+(2\ri+1)\lami\lamip^2\bigr)\nonumber\\
&\quad+n^{-2}\bigl(65\lami^2\lamip^2-(20\ri+10)\lami\lamip^2-(20\rip+50)\lami^2\lamip+(\rip^2+5\rip+6)\lami^2\nonumber\\
&\qquad\qquad+(\ri^2+\ri)\lamip^2+(4\ri\rip+10\ri+2\rip+5)\lami\lamip\bigr)\nonumber\\
&\quad-2\lami^2\lamip^2+n^{-1}\bigl(12\lami^2\lamip^2-(4\ri+2)\lami\lamip^2-(2\rip+4)\lami^2\lamip\bigr)\nonumber\\
&\quad+n^{-2}\bigl(-50\lami^2\lamip^2+(12\rip+24)\lami^2\lamip+(24\ri+12)\lami\lamip^2\nonumber\\
&\qquad\qquad-(4\ri\rip+8\ri+2\rip+4)\lami\lamip-(2\ri^2+2\ri)\lamip^2\bigr)\nonumber\\
&\quad+\lami^2\lamip^2+n^{-1}\bigl(-3\lami^2\lamip^2+(2\ri+1)\lami\lamip^2\bigr)\nonumber\\
&\quad+n^{-2}\bigl(7\lami^2\lamip^2-(6\ri+3)\lami\lamip^2+(\ri^2+\ri)\lamip^2\bigr)\nonumber\\
&\quad-2\lami^2\lamip^2+n^{-1}(12\lami^2\lamip^2-2\ri\lami\lamip^2-(4\rip+6)\lami^2\lamip\bigr)\nonumber\\
&\quad+n^{-2}\bigl(-50\lami^2\lamip^2+12\ri\lami\lamip^2+(24\rip+36)\lami^2\lamip\nonumber\\
&\qquad\qquad-(2\rip^2+6\rip+4)\lami^2-(4\ri\rip+6\ri)\lami\lamip\bigr)\nonumber\\
&\quad+4\lami^2\lamip^2+n^{-1}\bigl(-12\lami^2\lamip^2+4\ri\lami\lamip^2+(4\rip+4)\lami^2\lamip\bigr)\nonumber\\
&\quad+n^{-2}\bigl(-12\ri\lami\lamip^2-(12\rip+12)\lami^2\lamip+(4\ri\rip+4\ri)\lami\lamip+28\lami^2\lamip^2\bigr)\nonumber\\
&\quad-2\lami^2\lamip^2+n^{-1}\bigl(2\lami^2\lamip^2-2\ri\lami\lamip^2\bigr)+n^{-2}\bigl(2\ri\lami\lamip^2-2\lami^2\lamip^2\bigr)\nonumber\\
&\quad+\lami^2\lamip^2+n^{-1}\bigl(-3\lami^2\lamip^2+(2\rip+1)\lami^2\lamip\bigr)\nonumber\\
&\quad+n^{-2}\bigl((-6\rip-3)\lami^2\lamip+(\rip^2+\rip)\lami^2+7\lami^2\lamip^2\bigr)\nonumber\\
&\quad-2\lami^2\lamip^2+n^{-1}\bigl(2\lami^2\lamip^2-2\rip\lami^2\lamip\bigr)+n^{-2}\bigl(2\rip\lami^2\lamip-2\lami^2\lamip^2\bigr)\nonumber\\
&\quad+\lami^2\lamip^2+O(n^{-3})\nonumber\\
&=n^{-2}\bigl(3\lami^2\lamip^2-5\lami^2\lamip-\lami\lamip^2+2\lami^2+\lami\lamip\bigr)+O(n^{-3}).
\end{align}
Now we are ready to calculate \eqref{E_R2},\eqref{E_R3} and \eqref{E_R4}.
From \eqref{ex_form_E_R2} and 
\begin{align*}
&n^{-1}\bigl(E[\delip^2]+E[\deli^2]-2E[\deli\delip]\bigr)\\
&=n^{-1}\bigl(\lamip(1-\lamip)+\lami(1-\lami)-2\lami(1-\lamip)\bigr)\\
&\quad +n^{-2}\bigl(5\lamip^2-3\lamip-4\rip\lamip+\rip(\rip+1)+5\lami^2-3\lami-4\ri\lami+\ri(\ri+1)\\
&\qquad\qquad-10\lami\lamip+6\lami+4\rip\lami+4\ri\lamip-2\ri(1+\rip)\bigr)+O(n^{-3})\\
&=n^{-1}\bigl(-\lami+\lamip-(\lamip^2+\lami^2-2\lami\lamip)\bigr)\\
&\quad+n^{-2}\bigl(5(\lamip-\lami)^2-3(\lamip-\lami)-4(\rip-\ri)(\lamip-\lami)\\
&\qquad\qquad+\ri(\ri+1)+\rip(\rip+1)-2\ri(1+\rip)\bigr)+O(n^{-3})\\
&=n^{-1}(m_i-m_i^2)+n^{-2}(5m_i^2-3m_i-4(\rip-\ri)m_i+(\rip-\ri)(\rip-\ri+1)\bigr)+O(n^{-3}),
\end{align*}
we  have the result \eqref{E_R2} as follows;
\begin{align}
&\sum_{i=0}^p m_i E[R_i^2]\nonumber\\
&=\sum_{i=0}^{p}\bigl\{n^{-1}(1-m_i)+n^{-2}\bigl(5m_i-3-4(\rip-\ri)+(\rip-\ri)(\rip-\ri+1)m_i^{-1}\bigr)\bigr\}\nonumber\\
&\quad+O(n^{-3})\nonumber\\
&=n^{-1}\sum_{i=0}^p(1-m_i)+n^{-2}\bigl(\sum_{i=0}^p(5m_i-3)-4\sum_{i=0}^p(\rip-\ri)\nonumber\\
&\quad+\sum_{i=0}^p(\rip-\ri)(\rip-\ri+1)m_i^{-1}\bigr)+O(n^{-3})\nonumber\\
&=n^{-1}\sum_{i=0}^p(1-m_i)\nonumber\\
&\quad+n^{-2}\bigl(5-3(p+1)-4+\sum_{i=0}^{p}(\rip-\ri)(\rip-\ri+1)m_i^{-1}\bigr)+O(n^{-3})\nonumber\\
&=n^{-1}p+n^{-2}\bigl(-2-3p+\sum_{i=0}^{p}(\rip-\ri)(\rip-\ri+1)m_i^{-1}\bigr)+O(n^{-3}),
\end{align}
where we used the fact $\sum_{i=0}^p m_i=1$, $r_0=0$ and $r_{p+1}=1$.

From \eqref{ex_form_E_R3} and
\begin{align*}
&n^{-3/2}\bigl(E[\delip^3]-3E[\delip^2\deli]+3E[\delip\deli^2]-E[\deli^3]\bigr)\\
&=n^{-2}\bigl(7\lamip^3-3\lamip^2(\rip+3)+\lamip(3\rip+2)-21\lamip^2\lami+3\ri\lamip^2\nonumber\\
&\quad+(6\rip+27)\lami\lamip-3\ri\lamip-(6\rip+6)\lami+21\lami^2\lamip-(3\rip+18)\lami^2\nonumber\\
&\quad-(6\ri+9)\lami\lamip+(6\ri+3\rip+6)\lami-7\lami^3+3\lami^2(\ri+3)-\lami(3\ri,+2)\bigr)+O(n^{-3})\nonumber\\
&=n^{-2}\bigl[7(\lamip-\lami)^3+3(\lamip^2+\lami^2)(\ri-\rip-3)+6(\rip-\ri+3)\lami\lamip\nonumber\\
&\qquad\quad+\bigl(3(\rip-\ri)+2\bigr)(\lamip-\lami)\bigr]+O(n^{-3})\nonumber\\
&=n^{-2}\bigl[7(\lamip-\lami)^3+3(\ri-\rip-3)(\lamip-\lami)^2+\bigl(3(\rip-\ri)+2\bigr)(\lamip-\lami)\bigr]\nonumber\\
&\quad+O(n^{-3}),
\end{align*}
we have the result \eqref{E_R3} as follows;
\begin{align}
&\sum_{i=0}^p m_i E[R_i^3]\nonumber\\
&=n^{-2}\Bigl[\sum_{i=0}^{p}\bigl\{7m_i+3(\ri-\rip-3)+\bigl(3(\rip-\ri)+2\bigr)m_i^{-1}\bigr\}\Bigr]+O(n^{-3})\nonumber\\
&=n^{-2}\Bigl[\sum_{i=0}^p(7m_i-9)+3\sum_{i=0}^p(\ri-\rip)+\sum_{i=0}^p\bigl(3(\rip-\ri)+2)\bigr)m_i^{-1}\Bigr]+O(n^{-3})\nonumber\\
&=n^{-2}\Bigl[7-9(p+1)-3+\sum_{i=0}^p\bigl(3(\rip-\ri)+2)\bigr)m_i^{-1}\Bigr]+O(n^{-3})\nonumber\\
&=n^{-2}\Bigl[-5-9p+\sum_{i=0}^p\bigl(3(\rip-\ri)+2\bigr)m_i^{-1}\Bigr]+O(n^{-3}).
\end{align}

From \eqref{ex_form_E_R4} and
\begin{align*}
&n^{-2}\bigl(E[\delip^4]-4E[\delip^3\deli]+6E[\delip^2\deli^2]-4E[\deli^3\delip]+E[\deli^4]\bigr)\\
&=n^{-2}\bigl(3\lamip^4-6\lamip^3+3\lamip^2-12\lami\lamip^3+24\lamip^2\lami-12\lami\lamip\\
&\qquad\quad+18\lami^2\lamip^2-30\lami^2\lamip-6\lami\lamip^2+12\lami^2+6\lami\lamip-12\lami^3\lamip\\
&\qquad\quad+12\lami^3+12\lami^2\lamip-12\lami^2+3\lami^4-6\lami^3+3\lami^2\bigr)+O(n^{-3})\\
&=n^{-2}\bigl(3(\lamip^4-4\lami\lamip^3+6\lami^2\lamip^2-4\lami^3\lamip+\lami^4)\nonumber\\
&\qquad\quad-6(\lamip^3-3\lamip^2\lami+3\lami^2\lamip-\lami^3)+3(\lamip^2-2\lami\lami+\lami^2)\bigr)+O(n^{-3})\\
&=n^{-2}\bigl(3(\lamip-\lami)^4-6(\lamip-\lami)^3+3(\lamip-\lami)^2\bigr)+O(n^{-3}),
\end{align*}
we have the result \eqref{E_R4} as follows;
\begin{align}
&\sum_{i=0}^p m_i E[R_i^4]\nonumber\\
&=n^{-2}\Bigl[\sum_{i=0}^{p}(3m_i-6+3m_i^{-1})\Bigr]+O(n^{-3})\nonumber\\
&=n^{-2}\Bigl[3-6(p+1)+3\sum_{i=0}^p m_i^{-1}\Bigr]+O(n^{-3})\nonumber\\
&=n^{-2}\Bigl[-3-6p+3\sum_{i=0}^p m_i^{-1}\Bigr]+O(n^{-3}).
\end{align}
\\
\\
--\textit{Proof of \eqref{exp_ri}, \eqref{exp_ri^2}, \eqref{exp_ri_rip}}--
\\
\\
For $1\leq i \leq p$, 
\begin{align*}
E[\ri]&=\bar{r}_i(1+\bar{r}_i)+(1+\bar{r}_i)(-\bar{r}_i)=0,\\
E[\ri^2]&=\bar{r}_i^2(1+\bar{r}_i)-(1+\bar{r}_i)^2\bar{r}_i\\
&=(1+\bar{r}_i)\bar{r}_i\bigl(\bar{r}_i-(1+\bar{r}_i)\bigr)\\
&=-\bar{r}_i(1+\bar{r}_i),
\end{align*}
while $E[r_0]=0,\ E[r_{p+1}]=1, E[r_0^2]=0, \ E[r_{p+1}^2]=1$ is obvious from $r_i\equiv 0,\ r_{p+1}\equiv 1$.
\eqref{exp_ri_rip} is proved from \eqref{exp_ri} and the following equation; for $1\leq i \leq p-1$,
\begin{align*}
E[\ri(1-\rip)]&=\bari(1-\barip)(1+\bari)(1+\barip)+\bari(-\barip)(1+\bari)(-\barip)\\
&\quad+(1+\bari)(1-\barip)(-\bari)(1+\barip)+(1+\bari)(-\barip)(-\bari)(-\barip)\\
&=0.
\end{align*}

\end{document}